\documentclass[11pt]{amsart}
\usepackage{latexsym}
\usepackage{amssymb}
\usepackage{amscd}
\usepackage{mathrsfs}
\usepackage{rotating}
\newcommand\la{\lambda}
\newcommand\z{\zeta}
\newcommand\e{\eta}
\renewcommand\th{\theta}

\newcommand\io{\iota}
\newcommand\m{\mu}

\renewcommand\t{\tau}
\renewcommand\r{\rho}

\newcommand{\OO}{\mathbb O}

\newcommand\Ql{\bar{\mathbf Q}_l}

\newcommand\BP{\mathbf P}

\newcommand\BF{\mathbf F}
\newcommand\BC{\mathbf C}

\newcommand\BZ{\mathbf Z}

\newcommand\Bm{\mathbf m}

\newcommand\Bv{\mathbf v}
\newcommand\Bk{\mathbf k}

\newcommand\Bla{\boldsymbol\lambda}

\newcommand\Bmu{\boldsymbol\mu}

\newcommand\Bxi{\boldsymbol{\xi}}

\newcommand\CE{\mathcal{E}}

\newcommand\CP{\mathcal{P}}
\newcommand\CQ{\mathcal{Q}}

\newcommand\CX{ \mathcal{X}}

\newcommand\SB{\mathscr{B}}
\newcommand\SE{\mathscr{E}}

\newcommand\SL{\mathscr{L}}
\newcommand\SM{\mathscr{M}}
\newcommand\SN{\mathscr{N}} 
\newcommand\SO{\mathscr{O}}
\newcommand\SP{\mathscr{P}}
\newcommand\SQ{\mathscr{Q}}
\newcommand\SH{\mathscr{H}}

\newcommand\SX{\mathscr{X}}
\newcommand\SY{\mathscr{Y}}

\newcommand\Fg{\mathfrak g}

\newcommand\iv{^{-1}}
\newcommand\wh{\widehat}
\newcommand\wt{\widetilde}
\newcommand\wg{^{\wedge}}

\newcommand\ol{\overline}

\newcommand\lra{\leftrightarrow}

\newcommand\IC{\operatorname{IC}}

\newcommand\Hom{\operatorname{Hom}}
\newcommand\End{\operatorname{End}}

\newcommand\Ind{\operatorname{Ind}}

\newcommand\ch{\operatorname{ch}}

\newcommand\reg{_{\operatorname{reg}}}

\newcommand\unip{\operatorname{uni}}
\newcommand\uni{_{\operatorname{uni}}}
\newcommand\nil{_{\operatorname{nil}}}

\newcommand\lp{\operatorname{\!\langle\!}}
\newcommand\rp{\operatorname{\!\rangle\!}}

\newcommand\nat{^{\natural}}

\newcommand{\isom}{\,\raise2pt\hbox{$\underrightarrow{\sim}$}\,}
\numberwithin{equation}{section}

\newtheorem{thm}{Theorem}[section]
\newtheorem{lem}[thm]{Lemma}
\newtheorem{cor}[thm]{Corollary}
\newtheorem{prop}[thm]{Proposition}

\def \para#1{\par\medskip\textbf{#1}
              \addtocounter{thm}{1}}

\def \remark#1{\par\medskip\noindent
                \textbf{Remark #1}
                \addtocounter{thm}{1}}

\def \remarks#1{\par\medskip\noindent
                \textbf{Remarks #1}
                \addtocounter{thm}{1}}

\begin{document}
\setlength{\baselineskip}{4.0mm}
\setlength{\abovedisplayskip}{4.0mm}
\setlength{\belowdisplayskip}{4.0mm}
%%%
%%%
\renewcommand{\theenumi}{\roman{enumi}}
\renewcommand{\labelenumi}{(\theenumi)}
\renewcommand{\thefootnote}{\fnsymbol{footnote}}
%%%
\renewcommand{\thefootnote}{\fnsymbol{footnote}}
%%%
\allowdisplaybreaks[2]
%\NoBlackBoxes
\parindent=20pt
%\addtocounter{section}{1}

%%%%%%%%%%%%%%%%%%%%
%%%%%%%%%%%%%%%%%%%%%%%%%%%%%%%%%%%
\pagestyle{myheadings}
\medskip
\begin{center}
 {\bf Springer correspondence for complex reflection groups} 
\end{center}

\par\bigskip

\begin{center}
Toshiaki Shoji
\end{center}

\title{}

\begin{abstract}
This paper is a survey on the topics on the Springer 
correpsondence related to the varieties such as the enhanced 
variety or the exotic symmetric space. We explain in 
the case of exotic symmetric space of higher level, the 
complex reflection group $S_n\ltimes (\BZ/r\BZ)^n$ appears 
naturally in the framework of the Springer correspondence. 
\end{abstract}

\maketitle
%\markboth{SHOJI}{KOSTKA FUNCTIONS}
\pagestyle{myheadings}

\par\noindent
{\bf \S 1. Introduction}
\addtocounter{section}{1}

\par\medskip
Springer correspondence is a canonical correspondence between 
the  unipotent classes of a reductive group and the irreducible 
representations of its Weyl group, established by Springer [Sp]
in 1976.  In 1981, Lusztig [L1] found a way of reformulating Springer's 
theory in terms of the theory of perverse sheaves.  
In the same paper, he gave a geometric interpretation of Kostka 
polynomials in terms of the intersection cohomology associated to the 
unipotent classes in $GL_n$.  In 1980's, Lusztig established the theory 
of character sheaves, developing the idea in [L1].  
It became an interesting problem to generalize the theory of character 
sheaves to the case where the ambient variety is not a group.
Recently, by [AH], [FGT], [K1] and [SS1-3], interesting 
examples, such as the enhanced variety and the exotic symmetric space, 
were found, which enjoy a satisfied theory of character
sheaves and Springer correspondence, together with    
interesting relations with Kostka polynomials. 
\par
The exotic symmetric space and the enhanced variety of higher level 
are natural generalizations of those varieties mentioned above.  In [S3], 
the theory of Springer correspondence for the exotic symmetric space 
of level $r$ was established.  Interestingly, the complex reflection 
group $S_n\ltimes (\BZ/r\BZ)^n$ appears naturally in the description 
of the Springer correspondence. 
This paper is a survey on the Springer correspondence related to 
various varieties as above, based on the author's talk in the conference
of Tsinghua Sanya International Mathematics Forum in December 2014.   

\par\bigskip\bigskip
\noindent
{\bf S2. Springer correspondence for reductive groups}
\addtocounter{section}{1}

\para{2.1.}
First we review some hitorical results on the Springer correspondence.
The Springer correspondence is a natural correspondence between unipotent classes 
of a connected reductive group and irreducible reprensentations of the associated 
Weyl group, first estblished by Springer [Sp]. 
Here we give a formulation due to Lusztig [L1] and Borho-MacPherson [BM] based 
on the theory of perverse sheaves.
\par
Let $G$ be a connected reductive group over an algebraically closed field $\Bk$.
Let $B$ be a Borel subgroup of $G$ containing a maximal torus $T$ of $G$, 
$U$ the unipotent radical of $B$.  Let 
$W = N_G(T)/T$ be the Weyl group of $G$.  We denote by $\SB = G/B$ the flag variety 
of $G$. 
Consider the morphism 
\begin{equation*}
\pi : \wt G = \{ (x, gB) \in G \times \SB \mid g\iv xg \in B\} \to G
\end{equation*}
defined by $(x, gB) \mapsto x$. 
Then $\wt G$ is a smooth, irreducible variety, and $\pi$ is a proper map.
Let $G\reg$ be the set of regular semisimple elements in $G$, and put
$T\reg = T \cap G\reg$.       
Put $\wt G\reg = \pi\iv(G\reg)$ and let $\pi_0$ be the restriction of $\pi$ 
on $\wt G\reg$.  Here $T\reg \times G/T \isom \wt G\reg$  by 
$(x, gT) \mapsto (gxg\iv, gB)$, and $W$ acts on 
$\wt G\reg$ via the action $w : (x, gT) \mapsto (wxw\iv, gw\iv T)$ 
on $T\reg \times G/T$.   
Then $\pi_0 : \wt G\reg \to G\reg$ turns out to be a Galois covering with 
Galois group $W$. 
\par 
Consider a constant sheaf $\Ql$ on $G\reg$.  Then the direct image 
$(\pi_0)_*\Ql$ is a semisimple local system, such that 
$\End ((\pi_0)_*\Ql) \simeq \Ql[W]$, and is decomposed along 
the irreducible representations of $W$, 

\begin{equation*}
\tag{2.1.1}
(\pi_0)_*\Ql \simeq \bigoplus_{\r \in W\wg} \r \otimes \SL_{\r}, 
\end{equation*} 
where $W\wg$ denotes the set of irreducible representations of $W$ over $\Ql$, up to 
isomorphism, and 
$\SL_{\r} = \Hom_W ((\pi_0)_*\Ql, \r)$ is an irreducible local system on $G\reg$.
\par
Since $G\reg$ is an open dense smooth subset of $G$, one can consider 
the intersection cohomology complex $\IC(G, (\pi_0)_*\Ql)$ on $G$. 
Let $\pi_*\Ql$ be the direct image of the constant sheaf $\Ql$ on $\wt G$
(an object in the derived category of $\Ql$-sheaves).
It is proved by Lusztig [L1] that $\pi_*\Ql \simeq \IC(G, (\pi_0)_*\Ql)$. 
Hence $\pi_*\Ql[\dim G]$ is a semisimple perverse sheaf, equipped with $W$-action, 
and is decomposed as 

\begin{equation*}
\tag{2.1.2}
\pi_*\Ql[\dim G] \simeq \bigoplus_{\r \in W\wg} \r \otimes \IC(G, \SL_{\r})[\dim G].
\end{equation*}

\para{2.2.}
Let $G\uni$ be the set of unipotent elements in $G$, which is a closed subvariety of $G$, 
and is called the {\bf unipotent variety} of $G$. It is known that $G\uni$ consists of finitely 
many $G$-orbits, under the conjugation action of $G$. Let $\SN_G$ be the set of all the pairs 
$(C, \SE)$, where $C$ is a unipotent class, and $\SE$ is a $G$-equivariant simple local 
system on $C$.  If we fix $x \in C$, and let $A_G(x) = Z_G(x)/Z_G^0(x)$ be the component group 
of $Z_G(x)$.  Then the set of $G$-equivariant simple local system on $C$ is in bijection with 
the set $A_G(x)\wg$ of irreducible representations of (the finite group) $A_G(x)$. 
Thus $\SN_G$ can be written as

\begin{equation*}
\tag{2.2.1}
\SN_G \simeq \{ (x, \e) \mid x \in G\uni/\!\sim, \e \in A_G(x)\wg \}, 
\end{equation*} 
where $x$ runs over a set of representatives of $G$-orbits in $G\uni$.  
Let $\ol C$ be the closure of $C$ in $G\uni$, and consider the intersection 
cohomology $\IC(\ol C, \CE)$ on $\ol C$.  By the extension by zero, we regard 
$A_{(C,\SE)} = \IC(\ol C, \SE)[\dim C]$ as a perverse sheaf on $G\uni$.  Then   
$A_{(C,\SE)}$ is an $G$-equivariant simple perverse sheaf on $G\uni$, and 
the set of (isomorphism  class of) 
$G$-equivariant simple perverse sheaves on $G\uni$ is given by 
$\{ A_{(C,\SE)} \mid (C,\SE) \in \SN_G\}$. 
\par
Put $\wt G\uni = \pi\iv(G\uni)$, and let $\pi_1$ be the restriction of $\pi$ on $\wt G\uni$.
Thus 
\begin{equation*}
\pi_1 : \wt G\uni = \{ (x, gB) \in G\uni \times \SB \mid g\iv xg \in U\} \to G\uni
\end{equation*}  
with $\pi_1(x,gB) = x$. 
By the base change theorem, the restriction of $\pi_*\Ql$ on $G\uni$ is isomorphic to 
$(\pi_1)_*\Ql$, hence $(\pi_1)_*\Ql$ has a natural action of $W$. 
The following result is known as the Springer correspondence.

%%%%
%%%%
\begin{thm}[{Borho-MacPherson [BM]}]  %%%%  Theorem 2.3
Let the notations be as above. 
\begin{enumerate}
\item
$(\pi_1)_*\Ql[\dim G\uni]$ is a semisimple perverse sheaf om $G\uni$, 
equipped with $W$-action, and is decomposed as
\begin{equation*}
(\pi_1)_*\Ql[\dim G\uni] \simeq \bigoplus_{(C, \SE) \in \SN_G}
                    V_{(C,\SE)} \otimes A_{(C,\SE)},
\end{equation*} 
where $V_{(C,\SE)}$ is an irreducible representation of $W$ if it is non-zero.

\item
For any $\r \in W\wg$, there exists a unique pair $(C, \SE) \in \SN_G$ 
such that 

\begin{equation*}
\IC(G, \SL_{\r})|_{G\uni} \simeq \IC(\ol C, \SE)[\dim C - \dim G\uni]
\end{equation*}
and that $\r \simeq V_{(C, \SE)}$. 
\item 
The correspondence $\r \mapsto (C, \SE)$ in (ii) gives a bijection  
\begin{equation*}
W\wg \simeq \{(C,\SE) \in \SN_G \mid V_{(C,\SE)} \ne 0 \}. 
\end{equation*}
\end{enumerate}
\end{thm} 

\para{2.4.}
For $x \in G$, put $\SB_x = \{ gB \in \SB \mid g\iv xg \in B\}$.  Then $\SB_x$ 
is a closed subvariety of $\SB$.  Since $\pi\iv(x) \simeq \SB_x$, $\SB_x$ is called 
the {\bf Springer fibre} of $x$. Since $\pi_*\Ql$ is a complex with $W$-action, the stalk 
$\SH^i_x(\pi_*\Ql)$ at $x \in G$ of the cohomology sheaf $\SH^i(\pi_*\Ql)$ 
has a natural action of $W$. 
By a general theory $\SH^i_x(\pi_*\Ql) \simeq H^i(\SB_x, \Ql)$,  
so the cohomology group $H^i(\SB_x, \Ql)$ (a finite dimensional vector space over $\Ql$)  
has a structure of $W$-module.  The representation of $W$ on $H^i(\SB_x, \Ql)$ is 
called the {\bf Springer representation} of $W$. 
\par
Note that $Z_G(x)$ acts on $\SB_x$ by the left multiplication, and it induces an action 
of $Z_G(x)$ on $H^i(\SB_x, \Ql)$.  Since $Z_G(x)$ acts trivially on $H^i(\SB_x, \Ql)$, 
$A_G(x)$ acts on $H^i(\SB_x, \Ql)$.  
It is known that this action of $A_G(x)$ on $H^i(\SB_x, \Ql)$ commutes with the Springer 
action of $W$.  Thus $H^i(\SB_x,\Ql)$ has a structure of $W \times A_G(x)$-module.
Let $d_x = \dim \SB_x$.  Then it is known, for $x \in C$,  that 

\begin{equation*}
\tag{2.4.1}
d_x = \frac{1}{2}(\dim G\uni - \dim C).
\end{equation*}

\par
We now concentrate on the top cohomolgy group $H^{2d_x}(\SB_x, \Ql)$ as a $W \times A_G(x)$-module.
Theorem 2.3 can be reformulated in terms of the Springer representations of $W$,
which is the original form of the Springer correspondence due to Spriner [Sp]. 

\begin{cor}[{Springer [Sp]}]  %%%%%  Cor. 2.5
Write the decompostion of $W \times A_G(x)$-module $H^{2d_x}(\SB_x, \Ql)$ 
as follows;

\begin{equation*}
H^{2d_x}(\SB_x, \Ql) \simeq \bigoplus_{\e \in A_G(x)\wg}V_{(x,\e)} \otimes \e. 
\end{equation*}  
Then $V_{(x,\e)}$ is an irreducible $W$-module if it is non-zero.  Any irreducible 
representation of $W$ can be realized as $V_{(x, \e)}$ for a unique pair $(x, \e) \in \SN_G$
$($so $V_{(x,\e)}$ coincides with $V_{(C,\SE)}$ if $(x, \e) \lra (C,\SE)$ in (2.2.1) $)$.   
\end{cor}

\par\bigskip\noindent
{\bf \S 3. Kostka polynomials}

\addtocounter{section}{1}
\addtocounter{thm}{-5}

\para{3.1.}
We consider the Springer correspondence in Theorem 2.3 and Corollary 2.5 in the 
special case where $G = GL_n = GL(V)$. In this case, the set of unipotent classes in $G\uni$
is parametrized by the set $\SP_n$ of partitions of $n$, via Jorpdan normal form.
We denote by $C_{\la}$ the class corresponding to $\la \in \SP_n$.
In the case of $GL_n$, it is known that $Z_G(x)$ is connected, and so 
$A_G(x) = \{ 1\}$ for any $x \in G$.    
Thus $\SN_G  = \{ (C, \Ql) \mid \Ql: \text{ constant sheaf on } C \}$,
and $\SN_G$ can be identified with the set $G\uni/\!\sim$ of unipotent classes 
in $G$. 
Moreover, $W \simeq S_n$, the symmetric group of degree $n$, and the Springer 
correspondence gives a bijection $S_n\wg \simeq G\uni/\!\sim$. 
It is known that irreducible representations of $S_n$ are naturally parametrized 
by $\SP_n$.  We denote by $V_{\la}$ the irreducible representation of $S_n$ 
corresponding to $\la$ (here we use the labelling such that $V_{\la} = 1_{S_n}$ 
if $\la = (n)$, and $V_{\la}$ is the sign representation if $\la = (1^n)$.)  
Then it can be verified that the Springer correpsondence is actually given by 
$V_{\la} \lra C_{\la}$. 
In this case, the formula in Theorem 2.3 (i) is written as 

\begin{equation*}
\tag{3.1.1}
(\pi_1)_*\Ql[\dim G\uni] 
    \simeq \bigoplus_{\la \in \SP_n}V_{\la} \otimes \IC(\ol C_{\la}, \Ql)[\dim C_{\la}]
\end{equation*}

\par
Moreover, Corollary 2.5 implies,  
for $x \in C_{\la}$, that $H^{2d_x}(\SB_x, \Ql) \simeq V_{\la}$ as $S_n$-modules, and 
any irreducible representation of $S_n$ can be realized by such a Springer module. 

\par
Recall that the partial order, called the {\bf dominance order}, on $\SP_n$  
is defined as follows; 
for $\la, \mu \in \SP_n$, 
write  $\la = (\la_1, \dots, \la_k), \mu = (\mu_1, \dots, \mu_k) \in \SP_n$
for some common $k$ (by allowing  0 for the parts $\la_i, \mu_i$).  Then
$\mu \le \la$ if $\sum_{i = 1}^j\mu_i \le \sum_{i=1}^j\la_i$ for each $j$. 
\par
By using the dominance order, the closure $\ol C_{\la}$ can be described as fllows;

\begin{equation*}
\tag{3.1.2}
\ol C_{\la} = \coprod_{\mu \le \la}C_{\mu}.
\end{equation*}

\para{3.2.}
Kostka polynomials $K_{\la, \mu}(t) \in \BZ[t]$ are well-known polynomials 
in the combinatorial theory, indexed by partitions 
$\la,\mu \in \SP_n$.  They are given as the coefficients of the transition matrix 
between the basis of Schur functions and that of Hall-Littlewood functions in the 
space of symmetic functions.   
\par
For $\la = (\la_1, \dots \la_k) \in \SP_n$, we define an integer $n(\la) \in \BZ_{\ge 0}$ by 
$n(\la) = \sum_{i=1}^k(i-1)\la_i$. 
Then it is known that $K_{\la, \mu}(t)$ is a monic of degree $n(\mu) - n(\la)$ if $\mu \le \la$, 
and $K_{\la,\mu} = 0$ otherwise. 
We define a modified Kostka polynomial $\wt K_{\la, \mu}(t)$ by 
$\wt K_{\la,\mu}(t) = t^{n(\mu)}K_{\la, \mu}(t\iv)$. 
\par
The following interesting formula was proved by Lusztig in 1981.

%\addtocounter{thm}{-5}
\begin{thm}[Lusztig{[L1]}]
For $\la \in \SP_n$, consider the intersection cohomology $K = \IC(\ol C_{\la}, \Ql)$
on $G\uni$.  Let $\SH^i_xK$ be the stalk at $x \in G\uni$ of the cohomology shear $\SH^i K$.    
\begin{enumerate}
\item
$\SH^iK = 0$ for odd $i$. 

\item 
Assume that $\mu \le \la$, and take $x \in C_{\mu} \subset \ol C_{\la}$. 
Then

\begin{equation*}
\tag{3.3.1}
\wt K_{\la,\m}(t) = t^{n(\la)}\sum_{i \ge 0} (\dim \SH^{2i}_xK)t^i.
\end{equation*} 
\end{enumerate} 
\end{thm}

\remarks{3.4.} \ (i) Lascoux-Sch\"utzenberger theorem  ([M, III, (6.5)]) 
gives a combinatorial description 
of Kostka polynomials as follows;
\begin{equation*}
\tag{3.4.1}
K_{\la,\mu}(t) = \sum_{T \in SST(\la, \mu)}t^{c(T)},
\end{equation*}
where $SST(\la, \mu)$ is the set of semistandard tableaux of shape $\la$ and weight $\mu$, 
and $c(T)$ is the charge of the tableau $T$ (see [M] for the definition).
This thereom  implies that the coefficients of $K_{\la, \mu}(t)$
are non-negative integers.  Theorem 3.3 gives an alternate proof of this fact. 
\par
(ii) \ Assume that $\Bk$ is an algebraic closure of a finite field $\BF_q$, and 
consider a finite subgroup $GL_n(\BF_q)$ of $GL_n$.   
It was known that modified Kostka polynomials $\wt K_{\la, \mu}(t)$ 
evaluated at $t = q$ give some character values of certain irreducible representations 
of $GL_n(\BF_q)$ over $\Ql$.  
Then (3.3.1) is regarded as a formula which describes some character values of 
$GL_n(\BF_q)$ in terms of the intersection cohomology associated to $G$-orbits 
in $G$.  
This point of view was later generalized extensively by Lusztig, and he established 
the theory of character sheaves ([L3]) on connected reductive groups $G$, 
which is a geometric theory describing all the character values of $G(\BF_q)$ 
in terms of certain intersection cohomology associated to $G$.  
The theory of Springer correspondence was later generalized by Lusztig to  
the theory of generalized Springer correspondence ([L2]), which plays an essential 
role in the theory of character sheaves.

\par\bigskip\medskip
\noindent
{\bf \S 4.  Enhanced variety  $\SX^{\text{en}}$ }

\addtocounter{section}{1}
\addtocounter{thm}{-4}

\para{4.1.}
It is an interesting problem to generalize the theory of character sheaves 
to the case where the ambient variety is not a connected reductive group. 
In fact, in [L4], Lusztig developed the theory of character sheaves on 
disconnected reductive groups.    
Boyarchenko and Drinfeld [BD] developed the theroy of character sheaves on 
unipotent groups in positive characteristic. 
It is also interesting to replace $G$ by a variety $\SX$ on which $G$ acts. 
In fact, Ginzburg [Gi] defined the character sheaves on the symmetric space 
$G/H$ ($H$ is a closed subgroup of $G$). 
Recently, some other examples of $\SX$, such as the enhanced variety and the exotic 
symmetric space, were found that they enjoy a satisfied theory of character sheaves.

\para{4.2.}
Before going to the discussion on the enhanced variety and the exotic space, we 
prepare some general notion from the combinatorics. 
For a positive integer $r$, we denote by $\SP_{n,r}$ the set of $r$-tuple of 
partitions $\Bla = (\la^{(1)}, \dots, \la^{(r)})$ such that 
$\sum_{i=1}^r|\la^{(i)}| = n$. 
We shall define a partial order $\Bmu \le \Bla$ in $\SP_{n,r}$.
For $\Bla \in \SP_{n,r}$, express the partition $\la^{(i)}$ as 
$\la^{(i)} = (\la^{(i)}_1, \dots, \la^{(i)}_k)$ for some common $k$, 
by allowing zero on $\la^{(i)}_j$, and define a composition $c(\Bla)$ of $n$ 
by 

\begin{equation*}
c(\Bla) = (\la^{1)}_1, \dots, \la^{(r)}_1, \la^{(1)}_2,\dots, \la^{(r)}_2, 
\dots, \la^{(1)}_k, \dots, \la^{(r)}_k). 
\end{equation*} 
Then define $\Bmu \le \Bla$ by the condition  $c(\Bmu) \le c(\Bla)$, 
by using the dominance order on $\SP_n$ (can be defined similarly for 
the set of compositions of $n$).
\par
The $n$-function on $\SP_n$, $\la \mapsto n(\la)$, in 3.2 is also 
generalized to the case of $\SP_{n,r}$. We define a function 
$a: \SP_{n,r} \to \BZ_{\ge 0}$ as follows; for each $\Bla = (\la^{(1)}, \dots, \la^{(r)})$,  
put
\begin{equation*}
\tag{4.2.1}
a(\Bla) = r\cdot n(\Bla) + |\la^{(2)}| + 2|\la^{(3)}| + \cdots + (r-1)|\la^{(r)}|,
\end{equation*} 
where $n(\Bla) = \sum_{i=1}^r n(\la^{(i)})$. 
In particular, in the case where $r = 2$, $a(\Bla) = 2n(\Bla) + |\la^{(2)}|$. 
\para{4.3.}
First we shall introduce the enhanced variety. 
Let $G = GL_n = GL(V)$ be as in \S 3, where $V$ is an $n$-dimensiona vector 
space over $\Bk$.  We consider the direct product $\SX = G \times V$, on which 
$G$ acts by the conjugation action on the first factor, and the natural action on 
the second factor. The variety $\SX = \SX^{\text{en}}$ is called the 
{\bf enhanced variety}. 
Put $\SX\uni = G\uni \times V$, which is an $G$-stable closed subset of $G$. 
$\SX\uni$ has a role of the unipotent variety, and is called the unipotent variety 
of $\SX$. 
The geometry of the enhanced varieties $\SX$ and $\SX\uni$ is studied extensively 
by Achar-Henderson [AH], and Finkelberg-Ginzburg-Travkin [FGT].  In particular, 
the theory of character sheaves and the Springer correspondence on $\SX$ 
were discussed in [FGT].  
Note that $\SX\uni$ is isomorphic to the {\bf enhanced nilpotent cone} 
$\Fg\nil \times V$ (here $\Fg\nil$ is the nilpotent cone of the Lie algebra 
$\Fg$ of $G$) introduced by [AH].  
The following fact was found by 
 [AH] and Travkin [T], independently.

\begin{lem} [{[AH], [T]}] %%%  Lemma 4.4
Let $\SX\uni/\!\sim$ be the set of $G$-orbits in $\SX\uni$.  Then $\SX\uni$ is in bijection 
with $\SP_{n,2}$. 
\end{lem}      

\para{4.5.} Following [AH], we shall give an explicit correspondence. 
Take $z = (x,v) \in G\uni \times V$.  Put $E^x = \{ y \in \End(V) \mid xy = yx\}$.
$E^x$ is a subalgebra of $\End(V)$ containing $x$.  Then $V_z = E^xv$ is 
an $x$-stable subspace of $x$. Let $\la'$ (resp. $\la''$) be the Jordan type of 
$x|_{V_z}$ (resp. $x|_{V/V_z}$). 
We have $\Bla = (\la',\la'') \in \SP_{n,2}$, and the assignment $(x, v) \mapsto \Bla$ 
gives the required parametrization of $G$-orbits in $\SX\uni$.  We denote by 
$\SO_{\Bla}$ the $G$-orbit corresponding to $\Bla$. 
\par
The closure relations for $\ol\SO_{\Bla}$ were determined by [AH].

\begin{lem}[{[AH, Thm. 3.9]}]  %%%%  Lemma 4.6
For each $\Bla \in \SP_{n,2}$, 
\begin{equation*}
\ol\SO_{\Bla} = \coprod_{\Bmu \le \Bla}\SO_{\Bmu}.
\end{equation*}
\end{lem}

\para{4.7.}
Following [FGT], we shall describe the Springer correspondence for $\SX\uni$. 
We follow the notation in \S 3.  Let $(M_i)_{0 \le i \le n}$ be the total flag in $V$ whose 
stabilizer in $G$ equals to $B$.  Thus we can choose a basis $\{ e_1, \dots, e_n \}$ of $V$
consisitng of weight vectors for $T$ such that $M_i = \lp e_1, \dots, e_i \rp$.  
For an integer $m$ such that $0 \le m \le n$, we define 
\begin{align*}
\wt\SX_m &= \{ (x,v,gB) \in G \times V \times \SB \mid g\iv xg \in B, g\iv v \in M_m\}, \\
   \SX_m &= \bigcup_{g \in G}g(B \times M_m),
\end{align*} 
and define a map $\pi_m : \wt\SX_m \to \SX$ by $(x,v,gB) \mapsto (x,v)$.  
Then $\pi_m$ is a surjective map onto $\SX_m$.  Since $\pi_m$ is proper, 
$\SX_m$ is a closed subvariety of $\SX$. 
We also consider their restriction on the unipotent variety, 

\begin{align*}
\wt\SX_{m.\unip} &= \{(x,v, gB) \in G\uni \times V \times \SB 
             \mid g\iv xg \in U, g\iv v \in M_m \},  \\
\SX_{m,\unip} &= \bigcup_{g \in G}g(U \times M_m),  
\end{align*}
and define a map $\pi_{1,m}: \SX_{m,\unip} \to \SX_m$ by 
$(x,v,gB) \mapsto (x,v)$.  Then $\pi_{1,m}$ is a proper surjective map 
onto $\SX_{m,\unip}$. 
Put $\Bm = (m, n-m)$, and let $S_{\Bm} = S_m \times S_{n-m}$ be the Weyl subgroup of $S_n$, 
which is the stabilizer of $M_m$ in $S_n$ (under the embedding of $S_n$ into $G$ 
as a permutation group of $\{ e_1, \dots, e_n\}$). 
Put 
\begin{equation*}
M_m^0 = \{ v \in M_m \mid v = \sum_{i = 1}^ma_ie_i \text{ with } a_i \ne 0 
    \text{ for all } i \},
\end{equation*}
and put $\SY_m^0 = \bigcup_{g \in G}g(T\reg \times M_m^0)$.
Then $\SY^0_m$ is an open dense subset of $\SX_m$. Put $d_m = \dim \SX_m$.  
It is proved by Finkelberg and Ginzburg [FG, Cor. 5.4.2] that 
$(\pi_m)_*\Ql[d_m]$ is a semisimple perverse sheaf on $\SX_m$, 
equipped with $S_{\Bm}$-action, 
and is decomposed as 
\begin{equation*}
\tag{4.7.1}
(\pi_m)_*\Ql[d_m] \simeq \bigoplus_{\r \in S_{\Bm}\wg}\r \otimes \IC(\SX_m, \SL_{\r})[d_m],
\end{equation*}
where $\SL_{\r}$ is a simple local system  on $\SY^0_{m}$.  
\par
Let $\SP(\Bm)$ be the set of $\Bla = (\la^{(1)}, \la^{(2)}) \in \SP_{n,2}$ 
such that $|\la^{(1)}| = m, |\la^{(2)}| = n-m$. It is clear that irreducible representations of 
$S_{\Bm} = S_m \times S_{n-m}$ are parametrized by $\SP(\Bm)$.  We denote by 
$V_{\Bla}$ the irreducible representation of $S_{\Bm}$ corresponding to $\Bla$. 
\par
The following result was proved in [SS3, Thm. 2.11] 
(see also the proof of [FGT, Thm. 1]), 
which is regarded as an analogue of the 
Springer correspondence for the case $\SX\uni$.

\begin{thm} [{[SS3]}]  Put $d_m' = \dim \SX_{m, \unip}$. 
\begin{enumerate}
\item 
$(\pi_{1,m})_*\Ql[d'_m]$ is a semisimple perverse sheaf on $\SX_{m,\unip}$, 
equipped with $S_{\Bm}$-action, 
and is decomposed as 

\begin{equation*}
(\pi_{1,m})_*\Ql[d'_m] \simeq \bigoplus_{\Bla \in \SP(\Bm)}V_{\Bla} 
         \otimes \IC(\ol\SO_{\Bla}, \Ql)[\dim \SO_{\Bla}].
\end{equation*}

\item
For each $\Bla \in \SP(\Bm)$, let $\SL_{\r}$ be the 
simple local system on $\SY^0_m$ corresponding to $\r = V_{\Bla} \in S_{\Bm}\wg$.
Then we have

\begin{equation*}
\IC(\SX_m, \SL_{\r})|_{\SX_{m,\unip}}\simeq \IC(\ol\SO_{\Bla}, \Ql)[\dim \SO_{\Bla} - d'_m].
\end{equation*} 

\end{enumerate}
\end{thm}

\remarks{4.9.} 
(i) \ In this case, the Springer correspondence is given by the following diagram.

\begin{equation*}
\coprod_{0 \le m \le n}(S_m \times S_{n-m})\wg \simeq 
        \coprod_{0 \le m \le n}\{ \SO_{\Bla} \mid \Bla \in \SP(\Bm)\} = \SX\uni/\!\sim.
\end{equation*} 
\par
(ii) \ In the case of the enhanced variety, $Z_G(z)$ is connected for any $z \in \SX\uni$.
It follows that $A_G(z) = Z_G(z)/Z_G^0(z) = \{1\}$, and so the $G$-equivariant simple local system 
on the $G$-orbit is only the constant sheaf $\Ql$.  This situation is quite similar to the case
of $GL_n$ explained in \S 3.   

\para{4.10.}
In [S1,2], a generalization of Kostka polynomials was introduced, which are 
functions indexed by a pair of $r$-partitions of $n$.  They are apriori 
rational functions in $\BZ(t)$.  Here we consider such functions associated to 
``limit symbols'' as given in [S2, \S 3].
In the case where $r = 2$,  it was shown in [S2] that they are actually  
polynomials in $\BZ[t]$, which we denote by $K_{\Bla, \Bmu}(t)$ for 
$\Bla, \Bmu \in \SP_{n,2}$. (For the definition of $K_{\Bla, \Bmu}(t)$, see also [LS].)
By [S2, Prop. 3.3], $K_{\Bla, \Bmu}(t) = 0$ unless $\Bmu \le \Bla$, in which case
it is a monic of degree $a(\Bmu) - a(\Bla)$. So the 
modified Kostka polynomial is defined by 
$\wt K_{\Bla, \Bmu}(t) = t^{a(\Bmu)}K_{\Bla, \Bmu}(t\iv)$ as in the original case. 
The construction of $K_{\Bla, \Bmu}(t)$ is purely combinatorial, but Achar and Henderson 
proved in [AH, Thm. 5.2] that such Kostka polynomials can be interpreted 
by using the geometry of the 
enhanced variety $\SX\uni$, as in the case of $GL_n$.      
%%%
%%%
\begin{thm}[{Achar-Henderson [AH]}]  %%%%  Theorem 4.10
For $\Bla \in \SP_{n,2}$, consider the intersection cohomology $K = \IC(\ol\SO_{\Bla}, \Ql)$
on $\SX\uni$. 
\begin{enumerate}
\item
$\SH^i K = 0$ for odd $i$.
\item
Take $\Bla, \Bmu \in \SP_{n,2}$ such that $\Bmu \le \Bla$. 
Then for $z \in \SO_{\Bmu} \subset \ol\SO_{\Bla}$, we have
\begin{equation*}
\tag{4.10.1}
\wt K_{\Bla, \Bmu}(t) = t^{a(\Bla)}\sum_{i \ge 0}(\dim \SH^{2i}_zK)t^{2i}.
\end{equation*}
\end{enumerate}
\end{thm}

\par
Note that in the formula (4.10.1) each term $t^i$  
is replaced by $t^{2i}$, compared to the formula (3.3.1).

\par\bigskip\medskip
\noindent
{\bf \S 5.  Exotic symmetric space  $\SX^{\text{ex}}$ }

\addtocounter{section}{1}
\addtocounter{thm}{-11}

\para{5.1.}
Let $V$ be an $2n$-dimensional vector space over $\Bk$, where $\Bk$ is 
an algebraically closed field of $\ch \Bk \ne 2$. Put $G = GL_{2n} \simeq GL(V)$.
We consider the involutive automorphism $\th : G \to G$ defined by 
$\th(g) = J\iv({}^tg\iv)J$ with  
\begin{equation*}
J = \begin{pmatrix}
      0  &   1_n  \\
     -1_n  &  0   
    \end{pmatrix}. 
\end{equation*}
Let $H = G^{\th}$ be the group of $\th$-fixed elements in $G$.  Then 
$H$ coincides with the symplectic group $Sp_{2n}$. 
Let $\io : G \to G, g \mapsto g\iv$ be the anti-automorphim on $G$, and define 
a subset of $G$ by 
\begin{equation*}
G^{\io\th} = \{ g \in G \mid \th(g) = g\iv \}.
\end{equation*} 
$H$ acts on $G^{\io\th}$ by the conjugation action, and the map 
$g \to g\th(g)\iv$ induces an isomprphism $G/H \simeq G^{\io\th}$. 
Under this isomorphism, the left multiplication of $H$ on $G/H$ corresponds 
to the conjugation action of $H$ on $G^{\io\th}$. 
Thus, instead of considering the symmetric space $G/H$ with left $H$-action, 
we may consider the closed subvariety $G^{\io\th}$ of $G$ with conjugation 
action of $H$. 
Put $G^{\io\th}\uni = G^{\io\th} \cap G\uni$.  Then $G^{\io\th}\uni$ is an $H$-stable
subset of $G^{\io\th}$, which plays a role of the unipotent variety for $G/H$.
We consider the variety $\SX = G^{\io\th} \times V$ on which $H$ acts naturally, 
and put $\SX\uni = G^{\io\th}\uni \times V$. $\SX = \SX^{\text{ex}}$ 
is called the {\bf exotic symmetric space}, 
and $\SX\uni$ is an $H$-stable closed subset of $\SX$.
The geometry of $\SX$ and $\SX\uni$ was studied extensively by Kato [K1,2] and 
[SS1-3] from different points of view.  
In particular, the Springer correspondence was discussed in [K2] 
and [SS1]  The theory of character sheaves was developed in [SS1-3].  
Let $\Fg$ be the Lie algebra of $G$. $\th$ induces a linear involutive map 
$\th : \Fg \to \Fg$, and we denote by $\Fg^{-\th}$ the $-1$ eigenspace of $\Fg$.
Put $\Fg^{-\th}\nil = \Fg^{-\th} \cap \Fg\nil$, on which $H$ acts naturally.  
$\Fg^{-\th}\nil \times V$ 
is the {\bf exotic nilpotent cone} introduced by Kato [K1], and is isomorphic to 
$\SX\uni$ with $H$-action. 
The following lemma was proved in [K1].

\begin{lem}[{[K1]}]  %%%%  Lemma 5.2.
Let $\SX\uni/\!\sim$ be the set of $H$-orbits in $\SX\uni$.
Then $\SX\uni/\!\sim$ is in bijection with the set $\SP_{n,2}$. 
\end{lem}  

\para{5.3.}
 $G\uni \times V$ is an enhanced variety discussed in \S 4, and 
the set of $G$-orbits in $G\uni \times V$ is in bijection with 
$\SP_{2n,2}$.  $\SX\uni = G^{\io\th}\uni \times V$ is a subset of 
$\subset G\uni \times V$,  and the action of $H$ on $\SX\uni$ 
is compatible with the action of $G$ on $G\uni \times V$. 
The connection between $G$-orbits in $G\uni \times V$ and 
$H$-orbits in $\SX\uni$ was given by Achar-Henderson [AH, Thm. 6.1]
as follows; let $\OO_{\Bxi}$ be a $G$-orbit in $G\uni \times V$
with $\Bxi = (\xi',\xi'') \in \SP_{2n,2}$.  Then $\OO_{\Bxi} \cap \SX\uni \ne \emptyset$
if and only if $\Bxi$ is of the form  
$\xi' = \la'\cup \la', \xi'' = \la''\cup \la''$ for some 
$\Bla = (\la',\la'') \in \SP_n$, in which case $\OO_{\Bxi} \cap \SX\uni$ 
consists of a single $H$-orbit. We denote this $H$-orbit by $\SO_{\Bla}$.
Any $H$-orbit is obtained in this way, and this gives a parametrization of 
$H$-orbits in $\SX\uni$. 
Achar-Henderson also proved in [AH, Thm. 6.3], under this parametrization, 
that the closure relations for $\ol\SO_{\Bla}$ are given by a similar formula 
as Lemma 4.6, in terms of the partial order $\Bmu \le \Bla$ in $\SP_{n,2}$.   

\para{5.4.}
The Springer correspondence for $\SX\uni$ was first established by [K1,2] based on the 
Ginzburg theory on affine Hecke algebras. After that an alternate approach 
based on the theory of character sheaves 
was done by [SS1].  Here we follow the discussion 
in [SS1].
Let $T \subset B$ be a $\th$-stable pair of a maximal torus and a Borel subgroup of $G$.
Thus $T^{\th} \subset B^{\th}$ is a pair of a maximal torus and a Borel subgrou of $H$.   
We denote by $\SB^{\th} = H/B^{\th}$ the flag variety of $H$. 
We fix an isotropic flag $M_1 \subset M_2 \subset \cdots \subset M_n$ in $V$ 
stable by $B^{\th}$.
Define varieties

\begin{align*}
\wt\SX &= \{(x,v, gB^{\th}) \in G^{\io\th} \times V \times \SB^{\th}
             \mid g\iv xg \in B^{\io\th}, g\iv v \in M_n \}, \\
\wt\SX\uni &= \{ (x,v, gB^{\th}) \in \wt\SX \mid x \in G^{\io\th}\uni \}
\end{align*}
and define a map 
$\pi : \wt\SX \to \SX$ by $(x,v,gB^{\th}) \mapsto (x,v)$.  Then $\pi\iv(\SX\uni) = \wt\SX\uni$, 
and we define $\pi_1: \wt\SX\uni \to \SX\uni$ as the restriction of $\pi$ on $\wt\SX\uni$. 
$\pi,\pi_1$ are proper surjective maps, and $\wt\SX, \wt\SX\uni$ are smooth, irreducible 
varieties.  
\par
Let $W_n = S_n \ltimes (\BZ/2\BZ)^n$ be the Weyl group of type $C_n$.  Thus $W_n$ 
is the Weyl group of $H$. It is known that irreducible representations of $W_n$ are
(up to isomorphim) parametrized by $\SP_{n,2}$.  We denote by $\wt V_{\Bla}$ the irreducible 
representation of $W_n$ correpsonding to $\Bla \in \SP_{n,2}$. The following result was 
proved in [SS1].

\begin{thm}[{[SS1, Thm. 4.2]}]  %%%%   Theorem 5.5
$\pi_*\Ql[\dim \SX]$ is a semisimple perverse sheaf on $\SX$, 
equipped with $W_n$-action, and is decomposed as

\begin{equation*}
\tag{5.5.1}
\pi_*\Ql[\dim \SX] \simeq \bigoplus_{\Bla \in \SP_{n,2}}\wt V_{\Bla} \otimes A_{\Bla},
\end{equation*}
where $A_{\Bla}$ is an $H$-equivariant simple perverse sheaf on $\SX$. 
\end{thm}

\remarks{5.6.} (i) \ In the case of the enhanced variety,  for a fixed $m$, 
the local system $\SL_{\r}$ is obtained from the finite Galois covering 
$\pi_m\iv(\SY^0_m) \to \SY^0_m$ with Galois group $S_n \times S_{n-m}$.  Hence
in the decomposition (4.7.1), all the simple components 
$\IC(\SX_m, \SL_{\r})[\dim \SX_m]$ have the same support $\SX_m$.    
In the case of exotic symmetric space, for each $m$ such that 
$0 \le m \le n$, one can define a vareity $\SX_m$ by 
$\SX_m = \bigcup_{g \in H}g(B^{\th} \times M_m)$.  Then in the formula 
(5.5.1), the support of $A_{\Bla}$ runs over all $\SX_m$. 
Incidently, the construction of $W_n$-action on $\pi_*\Ql$ is more complicated 
than the enhanced case (see the discussion in 7.4). 
\par
(ii) Either in the enhanced case or the exotic case, we have a filtration 
$\SX_1 \subset \SX_2 \subset \cdots \subset \SX_n = \SX$. 
The map $\pi : \wt\SX \to \SX$ in the exotic case 
corresponds to the map $\pi_n$ in the enhanced case. 
In the exotic case, we have only to consider the map $\pi$, but in the 
enhaned case, the map $\pi_n$ does not reflect the general situation since the condition 
$g\iv v \in M_n = V$ is meaningless. 
\par\medskip
The following result was first proved by [K1, Thm. 8.3], [K2, Thm. G], and 
then reproved by [SS1, Thm. 5.4, Thm. 7.1] by a different method, 
which gives the Springer correspondence for the exotic case.

\begin{thm}[{[K1,2], [SS1]}]  %%%  Theorem 5.7
$(\pi_1)_*\Ql[\dim \SX\uni]$ is a semisimple perverse sheaf on $\SX\uni$, equipped 
with $W_n$-action, and is decomposed as 

\begin{equation*}
(\pi_1)_*\Ql[\dim \SX\uni] \simeq \bigoplus_{\Bla \in \SP_{n,2}}
        \wt V_{\Bla} \otimes \IC(\ol\SO_{\Bla},\Ql)[\dim \SO_{\Bla}].  
\end{equation*}
Moreover, $A_{\Bla}|_{\SX\uni} \simeq \IC(\ol\SO_{\Bla}, \Ql)$ up to shift. 
\end{thm}

\para{5.8.}  
We shall consider the Springer fibre for $\SX$.
For $z = (x,v) \in \SX$, put 
\begin{equation*}
\SB^{\th}_z = \{ gB^{\th} \in \SB^{\th} \mid g\iv xg \in B^{\io\th}, 
                             g\iv v \in M_n \}.
\end{equation*}
Then $\SB^{\th}_z$ is a closed subvariety of $\SB^{\th}$ isomorphic 
to $\pi\iv(z)$, and is called the Springer fibre of $z \in \SX$.
Put $d_{\Bla} = (\dim \SX\uni - \dim \SO_{\Bla})/2$. 
As a corollary to Theorem 5.7, we have

\begin{cor}
\begin{enumerate}
\item
$\dim \SB^{\th}_z = d_{\Bla}$ for $z \in \SO_{\Bla}$.

\item
$H^{2d_{\Bla}}(\SB^{\th}_z, \Ql) \simeq \wt V_{\Bla}$ as $W_n$-modules for $z \in \SO_{\Bla}$.
The assignment $\SO_{\Bla} \mapsto H^{2d_{\Bla}}(\SB^{\th}_z, \Ql)$ gives a bijective 
correpsondence
\begin{equation*}
\tag{5.9.1}
\SX\uni/\!\sim  \, \simeq W_n\wg.
\end{equation*}
\end{enumerate}
\end{cor}

\remark{5.10.} \ 
Compared to the enhanced case (see Remarks 4.9 (i)), 
the Springer correspondence (5.9.1) in the exotic case 
has a well-satisfied form.  This is because we can 
construct representations of $W_n$ in the exotic case, though
only representations of subgroups of $S_n$ in the enhanced case. 

\para{5.11.}
The relationship between the exotic symmetric space $\SX\uni$ and 
Kostka polynomials is studied in [K3] 
and [SS2].  Let $K_{\Bla, \Bmu}(t)$ be the Kostka polynomial 
associated to $\Bla, \Bmu \in \SP_{n,2}$ as discussed in 4.10.
We have the following result.

\begin{thm}[{[K3, Thm. E], [SS2, Thm. 5.7]}]  %%%%  Theorem 5.12.
For $\Bla \in \SP_{n,2}$, consider the intersection cohomology 
$K = \IC(\ol\SO_{\Bla}, \Ql)$ on $\SX\uni$.
\begin{enumerate}
\item
 $\SH^iK = 0$  unless $i \equiv 0 \pmod 4$.  
\item
Take $\Bla, \Bmu \in \SP_{n,2}$ such that $\Bmu \le \Bla$. 
Then for $z \in \SO_{\Bmu} \subset \ol\SO_{\Bla}$, we have

\begin{equation*}
\tag{5.12.1}
\wt K_{\Bla, \Bmu}(t) = t^{a(\Bla)}\sum_{i \ge 0}(\dim \SH^{2i}_zK)t^i.
\end{equation*}
\end{enumerate}
\end{thm}

\remarks{5.13.}
(i) \  The formula (5.12.1) was first conjectured by Achar and Henderson 
([AH, Conjecture 6.4]). 
An idea for the proof suggested by them  in [AH] was carried out 
by Kato (in the case where $\Bk = \BC$). 
He showed in [K3, Thm.~E], for each $z \in \SO_{\Bmu}$, 
 that the cohomology ring $H^{\bullet}(\SB^{\th}_z, \BC)$ has 
a De Concini-Procesi type interpretation as in the case of 
$GL_n$ ([DP]), i.e., there exists a graded algebra isomorphism 
between $H^{\bullet}(\SB^{\th}_z,\BC)$ and 
$R_{\Bmu} = \BC[x_1, \dots, x_n]/I_{\Bmu}$, 
compatible with the action of $W_n$, 
where $I_{\Bmu}$ is the ideal of all polynomials 
$p(x_1, \dots, x_n)$ such that 
$p(\partial/\partial x_1, \dots, \partial/\partial x_n)$
annihilates the Specht module $S_{\Bmu}$ realized in the 
homogeneous component of $\BC[x_1, \dots, x_n]$ of degree $a(\Bmu)$. 
Let $R^i_{\Bmu}$ be the $W_n$-module obtained as the $i$-th homogeneous 
part of $R_{\Bmu}$.  It was conjectured in [S2, 3.13] that 
for any irreducible representation $\wt V_{\Bla}$ of $W_n$, we have  
\begin{equation*}
\sum_{i \ge 0}\lp R^i_{\Bmu}, \wt V_{\Bla}\rp_{\,W_n}t^i = \wt K_{\Bla,\Bmu}(t).
\end{equation*}
Kato proved this conjecture in [K3, Thm. E], which  
provides a proof of (5.12.1), combined with 
the purity result ([K3, Cor. 5.3]). 
\par
The idea of the proof employed in [SS2], which was also suggested in [AH], 
is to construct an analogue of the theory of character sheaves on $\SX$, 
and to use the orthogonality relations for Green functions. 
\par
(ii) \ 
The discussion on $\SX\uni$ makes sense if we restrict ourselves to the symmetric space
$G^{\io\th}\uni \simeq \Fg^{-\th}\nil$ itself.    
In fact, the set of $H$-orbits in $G^{\io\th}$ is in bijection with the set $\SP_n$, 
and the orbit $\SO_{(-,\la'')}$ in $\SX\uni$ gives an $H$-orbit in $G^{\io\th}\uni$, which 
we denote by $\SO_{\la''}$ with $\la'' \in \SP_n$. 
Thus $K = \IC(\ol\SO_{\Bla}, \Ql)$ coincides with the intersection cohomology 
$K^{\text{sym}} = \IC(\ol\SO_{\la''}, \Ql)$ on $G^{\io\th}\uni$. 
Under this setup, the following result was proved by Henderson [H, Thm. 6,3], 
and reproved by [SS2, Thm. 5.10];  
$\SH^iK^{\text{sym}} = 0$ unless $i \equiv 0 \pmod 4$, and 
\begin{equation*}
\tag{5.13.1}
\wt K_{\la'', \mu''}(t^2) = t^{2n(\la'')}\sum_{i \ge 0}(\dim \SH^{2i}_xK^{\text{sym}})t^i
\end{equation*}
for $x \in \SO_{\mu''}$, 
where $\wt K_{\la'', \mu''}(t)$ is the original (modified) Kostka polynomial associated to 
$\la'', \mu'' \in \SP_n$.  
Note that the modulo 4 vanishing of the chomology sheaf $\SH^iK^{\text{sym}}$ 
was first noticed by Grojnowski in his thesis [Gr]. 
\par
Note that it is known by [AH, Cor. 5.3 (ii)] that 
$\wt K_{\Bla, \Bmu}(t) = t^n \wt K_{\la'', \mu''}(t^2)$ for such $\Bla, \Bmu$, 
hence (5.13.1) is obtained 
as the special case of (5.12.1).

\par\bigskip\medskip
\noindent
{\bf \S 6. Exotic symmetric space of higher level}
\addtocounter{section}{1}
\addtocounter{thm}{-13}

\para{6.1.}
We follow the notation in \S 5. 
For an integer $r \ge 2$, consider the varieties 

\begin{equation*}
\SX' = G^{\io\th} \times V^{r-1} \supset
\SX'\uni = G^{\io\th}\uni \times V^{r-1}
\end{equation*}
with diagonal action of $H$ on $\SX', \SX'\uni$.
$\SX', \SX'\uni$ are a natural generalization of the exotic 
symmetric space studied in \S 5. 
But it occurs a crucial difference when  we consider the general $r$, i.e., 
\par\medskip
\noindent
(6.1.1) \ If $r \ge 3$, $\SX'\uni$ has infinitely many $H$-orbits.

\par\medskip
In fact, since 
$\dim G^{\io\th}\uni = 2n^2 - 2n$, 
we have

\begin{equation*}
\dim \SX'\uni = 2n^2 -2n + (r-1)2n > \dim H = 2n^2 + n 
\end{equation*}
if $r \ge 3$.
\par
Let us define varieties 

\begin{align*}
\wt\SX  &= \{ (x, \Bv, gB^{\th}) \in G^{\io\th} \times V^{r-1} \times \SB^{\th}
                \mid g\iv xg \in B^{\io\th}, g\iv\Bv \in M_n^{r-1} \},  \\
\SX &= \bigcup_{g \in H}g(B^{\io\th} \times M_n^{r-1}),  \\
\wt\SX\uni &= \{ (x, \Bv, gB^{\th}) \in \wt\SX \mid x \in G^{\io\th}\uni \}, \\
\SX\uni &= \{(x, \Bv) \in \SX \mid x \in G^{\io\th}\uni \} 
\end{align*}
and a map $\pi : \wt\SX \to \SX$ by $(x, \Bv, gB^{\th}) \mapsto (x,\Bv)$.
We also define $\pi_1 : \wt\SX\uni \to \SX\uni$ by the restriction of 
$\pi$ on $\wt G\uni$. 
The maps $\pi, \pi_1$ are proper surjective.  Hence $\SX$ is a closed subset of 
$\SX' = G^{\io\th} \times V^{r-1}$, and simiarly for $\SX\uni$.
$\SX$ is called the {\bf exotic symmetric space of level $r$}.
Note that the map $\pi : \wt\SX \to \SX'$ is not necessarily surjective if $r \ge 3$.
So we replace $\SX'$ by the image $\pi(\wt\SX) = \SX$, and $\SX'\uni$ by 
$\pi_1(\wt\SX\uni)$.  But $\SX\uni$ has still infinitely many $H$-orbits.   
The Springer correspondence for $\SX\uni$ was established in [S3], which will be 
discussed in next section.  In this section, we prepare some notations.  

\para{6.2.}
We consider certain subvarieties of $\SX, \SX\uni$, as in the case of the enhanced variety.
Put 
\begin{align*}
\SQ_{n,r} &= \{ \Bm = (m_1, \dots, m_r) \in \BZ_{\ge 0}^r \mid \sum m_i = n\}, \\
\SQ_{n,r}^0 &= \{ \Bm \in \SQ_{n,r} \mid m_r = 0 \}.
\end{align*}
By fixing $\Bm \in \CQ_{n,r}$, define $p_1, \dots, p_r$ by 
$p_k = m_1 + \cdots + m_k$. 

\par
For each $\Bm \in \CQ_{n,r}$, put
\begin{align*}
\wt\SX_{\Bm} &= \{ (x, \Bv, gB^{\th}) \in G^{\io\th} \times V^{r-1} \times \SB^{\th} 
                     \mid g\iv xg \in B^{\io\th}, g\iv \Bv \in \prod_{i=1}^{r-1}M_{p_i} \} \\
\SX_{\Bm} &= \bigcup_{g \in H}g(B^{\io\th} \times \prod_{i=1}^{r-1}M_{p_i}) \\
\pi^{(\Bm)} &: \wt\SX_{\Bm} \to \SX_{\Bm} \quad (x, \Bv, gB^{\th}) \mapsto (x, \Bv)  
\end{align*}
\par
$\wt\SX_{\Bm}$ is smooth, irreducible, and $\pi^{(\Bm)}$ is proper, surjecitve. 
Hence $\SX_{\Bm}$ is a closed subset of $G^{\io\th} \times V^{r-1}$.
Note that if
$\Bm = (n, 0, \dots, 0)$, then $\wt\SX_{\Bm} = \wt\SX$ and $\SX_{\Bm} = \SX$. 
The dimension of the varieties $\wt\SX_m, \SX_m$ are computed as follows;

\begin{lem} %%%  Lemma 6.2  
\begin{enumerate}
\item
$\dim \wt\SX_{\Bm} = 2n^2 + \sum_{i=1}^r(r-i)m_i$.

\item
$\dim \SX_{\Bm} = 2n^2 + \sum_{i = 1}^r (r-i)m_i - m_r$.
In particilar, if $m_r = 0$, i.e., if $\Bm \in \SQ^0_{n,r}$, 
then $\dim \wt\SX_{\Bm} = \dim \SX_{\Bm}$. 
\end{enumerate}
\end{lem}
\par
The condition $\dim \wt\SX_{\Bm} = \dim\SX_{\Bm}$ plays an important role 
in later discussions (the semi-smallness of the map $\pi^{(\Bm)}$). 
So in order to guarantee this condition, we pose the assumption 
$\Bm \in \SQ^0_{n,r}$ in some situations. 

\para{6.4.}
Let $W_{n,r} = S_n \ltimes (\BZ/r\BZ)^n$ be the complex reflection group 
$G(r,1,n)$. 
It is well-known that the set $W_{n,r}\wg$ of irreducible representations of
$W_{n,r}$ is parametrized by the set $\SP_{n,r}$. We denote by $\wt V_{\Bla}$ 
the irreducible representation of $W_{n,r}$ corresponding to $\Bla \in \SP_{n,r}$. 
For later use, we review the construction of $\wt V_{\Bla}$.  
Let $\z \in \Ql^*$ be a primitive $r$-th root of unity.   
Define a linear character $\t_i : \BZ/r\BZ \to \Ql^*$ by $1 + r\BZ \mapsto \z^{i-1}$  
for $1 \le i \le r$. 
For $\Bm \in \CQ_{n,r}$, put 
$S_{\Bm} = S_{m_1} \times S_{m_2} \times \cdots \times S_{m_r} \subset S_n$.
For each $i$, put $\wt S_{m_i} = S_{m_i} \ltimes (\BZ/r\BZ)^{m_i}$, and consider 
a subgroup 
$\wt S_{\Bm} = \wt S_{m_1} \times \cdots \times \wt S_{m_r}$
of $W_{n,r}$.
\par
An irreducible representation $\r$ of $S_{\Bm}$ can be written as 
$\r = \r_1\boxtimes\cdots \boxtimes \r_r$ with $\r_i \in S_{m_i}\wg$. 
We extend $\r_i$ to the irreducible representation $\wt\r_i$ on $\wt S_{m_i}$
by defining the action of $(\BZ/r\BZ)^{m_i}$ on the space 
$\r_i$ by $\t_i^{\otimes m_i}$.  
Let  
$\wt\r = \wt\r_1 \boxtimes \cdots \boxtimes \wt\r_r \in \wt S_{\Bm}\wg$,
and put 
\begin{equation*}
\wt V_{\r} = \Ind_{\wt S_{\Bm}}^{W_{n,r}}\wt\r. 
\end{equation*}
Then $\wt V_{\r}$ gives an irreducible representation of $W_{n,r}$. 
Since $S_{m_i}\wg \simeq \SP_{m_i}$, we can write $\r_i = V_{\la^{(i)}}$ 
for $\la^{(i)} \in \SP_{m_i}$. 
Then $\Bla = (\la^{(1)}, \dots, \la^{(r)}) \in \SP_{n,r}$, and we write 
$\r = \r_{\Bla}$ and $\wt V_{\r} = \wt V_{\Bla}$.  
This gives the required parametriation $W_{n,r}\wg \simeq \SP_{n,r}$. 

\para{6.5.}
Here we introduce a partition of $W_{n,r}\wg$ corresponding to $\SQ^0_{n,r}$.
Assume that $\Bm \in \SQ^0_{n,r}$, namely $\Bm = (m_1, \dots, m_{r-1}, 0)$.
For $0 \le k \le m_{r-1}$, put $\Bm(k) = (m_1, \dots, m_{r-2}, k, k') \in \SQ_{n,r}$
with $k + k' = m_{r-1}$.
Put 
$S_{\Bm(k)} = S_{m_1} \times \cdots \times S_{m_{r-2}} \times S_k \times S_{k'}$,
and define 
\begin{equation*}
(W_{n,r}\wg)_{\Bm} = \{ \wt V_{\r} \mid \r \in S\wg_{\Bm(k)} \text{ for } 0 \le k \le m_{r-1} \}
\subset W_{n,r}\wg.
\end{equation*}
Then we have a partition of $W_{n,r}\wg$ 

\begin{equation*}
W_{n,r}\wg = \coprod_{\Bm \in \SQ^0_{n,r}}(W\wg_{n,r})_{\Bm}.
\end{equation*}

For $\Bm \in \SQ_{n,r}$, let 
$\SP(\Bm)$ be the subset of $\SP_{n,r}$ consisting of $\Bla$ such that 
$|\la^{(i)}| = m_i$ for $i = 1, \dots, r$. 
For $\Bm \in \SQ_{n,r}^0$, put

\begin{equation*}
\wt\SP(\Bm) = \coprod_{0 \le k \le m_{r-1}}\SP(\Bm(k)).
\end{equation*}
Thus $\wt\SP(\Bm)$ is the subset of $\SP_{n,r}$ consisting of 
$\Bla$ such that $|\la^{(i)}| = i$ for $i = 1, \dots, r-2$. 
It is easy to see that 
\begin{equation*}
\tag{6.5.1}
(W_{n,r}\wg)_{\Bm} = \{ \wt V_{\Bla} \mid \Bla \in \wt\SP(\Bm)\}.
\end{equation*}

\para{6.6.}
Recall that $W_n$ ($= W_{n,2}$ in the notation of 6.4) is the Weyl group of type $C_n$.
For $\Bm \in \SQ^0_{n,r}$, we define a parabolic subgroup $W_{\Bm}\nat$ of $W_n$ by

\begin{equation*}
W_{\Bm}\nat = S_{m_1} \times \cdots \times S_{m_{r-2}} \times W_{m_{r-1}}  \subset W_n.
\end{equation*} 
For $\r = \r_1 \boxtimes \cdots \boxtimes \r_r \in S_{\Bm(k)}\wg$
we define an irreducible $W_{\Bm}\nat$-module $V_{\r}\nat$ by 

\begin{equation*}
V_{\r}\nat = \r_1 \boxtimes \cdots \boxtimes \r_{r-2} \boxtimes \psi_{r-1},
\end{equation*}
where 
\begin{equation*}
\psi_{r-1} = \Ind_{\wt S_k \times \wt S_{k'}}^{W_{m_{r-1}}}
         (\wt\r_{r-1} \boxtimes \wt\r_r)
\end{equation*}
(apply the construction of $\wt V_{\Bla}$ for the case $r = 2$).
Then we have a natural bijection 
\begin{equation*}
(W_{\Bm}\nat)\wg \simeq \coprod_{0 \le k \le m_{r-1}}S_{\Bm(k)}\wg \simeq (W_{n,r}\wg)_{\Bm}
\end{equation*}
through 
$V_{\r}\nat \lra \r  \lra \wt V_{\r}$.
We denote by $V\nat_{\Bla}$ the irreducible representation of $W\nat_{n,r}$ 
corresponding to $\wt V_{\Bla}$.  It is easy to see that 

\begin{equation*}
\tag{6.6.1}
(W\nat_{\Bm})\wg  = \{ V\nat_{\Bla} \mid \Bla \in \wt\SP(\Bm) \}.
\end{equation*} 

\bigskip\medskip
\noindent
{\bf \S 7. Springer correspondence for $\SX^{\text{ex}}$ of higher level}

\addtocounter{section}{1}
\addtocounter{thm}{-6}

\para{7.1}
In this section, we shall discuss about the Springer correspondence for 
the exotic symmetric space $\SX$ of level $r$ based on [S3]. 
First we generalize Theorem 5.5 to the case where $r$ is arbitary. 
Take $\Bm \in \SQ^0_{n,r}$.  
Recall the map 
$\pi^{(\Bm)} : \wt\SX_{\Bm} \to \SX_{\Bm}$. 
We define a map  $\pi_{\Bm} : \pi\iv(\SX_{\Bm}) \to \SX_{\Bm}$ as the 
restriction of $\pi$ on $\pi\iv(\SX_{\Bm})$.
Note that $\wt\SX_{\Bm} \subset \pi\iv(\SX_{\Bm})$ since 

\begin{align*}
\pi\iv(\SX_{\Bm}) &= \{ (x,\Bv, gB^{\th})\mid  (x, \Bv) \in \SX_{\Bm}, 
         g\iv xg \in B^{\io\th}, g\iv \Bv \in M_n^{r-1}\}, \\
\wt\SX_{\Bm} &= \{ (x, \Bv, gB^{\th}) \in \pi\iv(\SX_{\Bm}) \mid 
            g\iv\Bv \in \prod_{i=1}^{r-1}M_{p_i} \}.
\end{align*}

Put $d_{\Bm} = \dim \SX_{\Bm}$.  The following two results are both 
generalizations of Theorem 5.5 in the case $r = 2$.
(Note that $\pi_{\Bm}$ in this paper is written as $\ol\pi_{\Bm}$ in [S3]). 

\begin{thm}[{[S3, Thm. 3.2]}]   %%%%  Theorem 7.2
Assume that $\Bm \in \SQ^0_{n,r}$.  Then $\pi^{(\Bm)}_*\Ql[d_{\Bm}]$ is a semisimple perverse sheaf
on $\SX_{\Bm}$, 
equipped with $W_{\Bm}\nat$-action, and is decomposed as

\begin{equation*}
\pi^{(\Bm)}_*\Ql[d_{\Bm}] \simeq \bigoplus_{0 \le k \le m_{r-1}}
           \bigoplus_{\r \in S_{\Bm(k)}\wg} V_{\r}\nat \otimes 
             \IC(\SX_{\Bm(k)}, \SL_{\r})[d_{\Bm(k)}],
\end{equation*}
where $\SL_{\r}$ is a simple local system on a ceertain open dense subset of $\SX_{\Bm(k)}$
associated to $\r \in S_{\Bm(k)}\wg$.  
\end{thm}

\begin{thm}[{[S3, Thm. 2.2]}]  %%%%  Theorem 7.3. 
Assume that $\Bm \in \SQ^0_{n,r}$.  Then $(\pi_{\Bm})_*\Ql[d_{\Bm}]$ is a semisimple perverse 
sheaf on $\SX_{\Bm}$, equipped with $W_{n,r}$-action, and is decomposed as 

\begin{equation*}
(\pi_{\Bm})_*\Ql[d_{\Bm}] \simeq \bigoplus_{0 \le k \le m_{r-1}}
           \bigoplus_{\r \in S_{\Bm(k)}\wg}\wt V_{\r} \otimes 
             \IC(\SX_{\Bm(k)}, \SL_{\r})[d_{\Bm(k)}].
\end{equation*}
\end{thm}

\para{7.4.}
The construction of $W_{n,r}$-action on $(\pi_{\Bm})_*\Ql$ is a natural 
generalization of the construction of $W_n$-action in the case where $r = 2$.
Here we give 
some explanation on the construction of $W_{n,r}$-action and 
on the definition of local systems $\SL_{\r}$ involved in the theorems.
\par
Let $T^{\io\th}\reg$ be the set of regular semisimple elements in $T^{\io\th}$, 
i.e., the set of elements such that all the eigenvalues have multiplicity 2.  
Let 
\begin{equation*}
G^{\io\th}\reg = \bigcup_{g \in H}gT^{\io\th}\reg g\iv \subset G^{\io\th}
\end{equation*}
be the set of regular semisimple elements in $G^{\io\th}$, which is an open 
dense subset of $G^{\io\th}$. For $\Bm \in \SQ_{n,r}$, we define 
\begin{align*}
\wt\SY_{\Bm} &= \{ (x, \Bv, gB^{\th}) \in \wt\SX_{\Bm} \mid x \in G^{\io\th}\reg \}, \\
\SY_{\Bm} &=  \{ (x, \Bv) \in \SX_{\Bm} \mid x \in G^{\io\th}\reg \},  \\
\psi^{(\Bm)} &: \wt\SY_{\Bm} \to \SY_{\Bm} \quad (x, \Bv, gB^{\th}) \mapsto (x, \Bv).
\end{align*}
Then $\SY_{\Bm}$ is an open dense subset in $\SX_{\Bm}$, and  
$\psi^{(\Bm)}$ is the restriction of $\pi^{(\Bm)}$ on $\wt\SY_{\Bm}$.
We put $\psi : \wt \SY \to \SY$, where $\wt\SY = \wt\SY_{\Bm}, 
\SY = \SY_{\Bm}$ for $\Bm = (n, 0, \dots, 0)$.  
As in the case of $\pi_{\Bm}$, we define  
$\psi_{\Bm} : \psi\iv(\SY_{\Bm}) \to \SY_{\Bm}$ as the restriction of $\psi$ on 
$\psi\iv(\SY_{\Bm})$. Hence $\wt\SY_{\Bm} \subset \psi\iv(\SY_{\Bm})$.
Note that 
$\wt\SY_{\Bm}$ is expressed as 
\begin{equation*}
\wt\SY_{\Bm} \simeq H \times^{B^{\th} \cap Z_H(T^{\io\th})}
              (T^{\io\th}\reg \times \prod_i M_{p_i}).
\end{equation*}

Define a partial order on $\SQ_{n,r}$ by 
$\Bm' \le \Bm$ if and only if $p_i' \le p_i \  (1 \le i \le r)$,
where $p_i = m_1 + \cdots + m_i, p_i' = m_1' + \cdots + m_i'$. 
Define $\SY_{\Bm}^0 = \SY_{\Bm} \ \backslash \ \bigcup_{\Bm' < \Bm}\SY_{\Bm'}$. 
Then $\SY_{\Bm}^0$ is an open dense subset of $\SY_{\Bm}$. 
We define a subvariety $\wt\SY_{\Bm}^0$ of $\wt\SY_{\Bm}$ by 

\begin{equation*}
\wt\SY_{\Bm}^0 = H \times^{B^{\th} \cap Z_H(T^{\io\th})}
                      (T^{\io\th}\reg \times \prod_iM_{p_i}^0),
\end{equation*}
where $M_{p_i}^0$ is defined as in 4.7.  
$S_n$ acts naturally on $\psi\iv(\SY^0_{\Bm})$, and 
we have 
\begin{equation*}
\tag{7.4.1}
\psi\iv(\SY_{\Bm}^0) \simeq \coprod_{w \in S_n/S_{\Bm}}w(\wt\SY_{\Bm}^0),
\end{equation*}
which gives the decomposition of $\psi\iv(\SY^0_{\Bm})$ into irreducible components. 
\par
We define a variety 
\begin{equation*}
\wh\SY^0_{\Bm} = H \times^{Z_H(T^{\io\th})_{\Bm}}(T^{\io\th}\reg \times \prod_iM^0_{p_i}), 
\end{equation*}
where
\begin{align*}
Z_H(T^{\io\th}) &\simeq SL_2 \times \cdots \times SL_2  \quad \text{ ($n$-times) },   \\
Z_H(T^{\io\th}) \cap B^{\th} &\simeq  B_2 \times \cdots \times B_2 
                    \quad \text{ ($B_2$: Borel of $SL_2$) },  \\
Z_H(T^{\io\th})_{\Bm} &\simeq B_2 \times \cdots B_2 \times 
       \underbrace{SL_2 \times \cdots \times SL_2}_{\text{last $m_r$-factors}}.  
\end{align*}
Then the map $\psi^{(\Bm)} : \wt\SY^0_{\Bm} \to \SY^0_{\Bm}$ is decomposed as 
\begin{equation*}
\begin{CD}
\psi^{(\Bm)} : \wt\SY^0_{\Bm} @>\xi_{\Bm}>> \wh\SY_{\Bm}^0  @>\e_{\Bm}>>  \SY^0_{\Bm},
\end{CD}
\end{equation*}
where 
$\xi_{\Bm}$ is a locally trivial fibration with fibre $\simeq \BP_1^{m_r}$, and
$\e_{\Bm}$ is a finite Galois covering with group $S_{\Bm}$. 
It follows that 

\begin{align*}
\tag{7.4.2}
(\xi_{\Bm})_*\Ql &\simeq H^{\bullet}(\BP_1^{m_r})\otimes \Ql, \\  
\tag{7.4.3}
(\e_{\Bm})_*\Ql &\simeq \bigoplus_{\r \in S_{\Bm}\wg}\r \otimes \SL_{\r}, 
\end{align*}
where $\SL_{\r}$ is a simple local sytem on $\SY^0_{\Bm}$. 
We define $\psi_{\Bm}^0 : \psi\iv(\SY^0_{\Bm}) \to \SY_{\Bm}^0$ by the restriciton of 
$\psi : \wt\SY \to \SY$. 
In view of (7.4.1) $\sim$ (7.4.3), we have

\begin{align*}
(\psi^0_{\Bm})_*\Ql &\simeq \bigoplus_{\r \in S_{\Bm}\wg}\Ind_{S_{\Bm}}^{S_n}
                  (H^{\bullet}(\BP_1^{m_r}) \otimes \r) \otimes \SL_{\r}   \\
                    &\simeq \bigoplus_{\r \in S_{\Bm}\wg}\Ind_{\wt S_{\Bm}}^{W_{n,r}}
                  (H^{\bullet}(\BP_1^{m_r}) \otimes \wt\r') \otimes \SL_{\r}.
\end{align*}
In the second formula, one can define an  action of $(\BZ/r\BZ)^n$ on 
$H^{\bullet}(\BP_1^{m_r})$ and on $\r$ so that 
$H^{\text{top}}(\BP_1^{m_r})\otimes \wt \r' \simeq \wt\r \in \wt S_{\Bm}\wg$.
Thus $(\psi^0_{\Bm})_*\Ql$ turns out to be a local system equipped with $W_{n,r}$-action.
By using the decomposition $\SY_{\Bm} = \coprod_{\Bm' \le \Bm}\SY^0_{\Bm'}$, we can 
determine the decomposition of $(\psi_{\Bm})_*\Ql[d_{\Bm}]$ as follows;

\begin{equation*}
\tag{7.4.4}
(\psi_{\Bm})_*\Ql[d_{\Bm}] \simeq \bigoplus_{0 \le k \le m_{r-1}}\bigoplus_{\r \in S_{\Bm(k)}\wg}
       \wt V_{\r} \otimes \IC(\SY_{\Bm(k)}, \SL_{\r})[d_{\Bm(k)}].
\end{equation*}
This formula corresponds to the formula (2.1.1) in the original case, and one 
can prove that $(\pi_{\Bm})_*\Ql[d_{\Bm}]$ coincides with the  
intermidate extension of $(\psi_{\Bm})_*\Ql[d_{\Bm}]$ (note that $\SY_{\Bm(k)}$ 
is open dense in $\SX_{\Bm(k)}$).  

\para{7.5.}
We now discuss about the unipotent variety $\SX\uni$. 
Recall the map $\pi_1 : \wt\SX\uni \to \SX\uni$.
For $\Bm \in \SQ_{n,r}$, put

\begin{align*}
\wt\SX_{\Bm,\unip} &= \{ (x, \Bv, gB^{\th}) \in \wt\SX_{\Bm} \mid x \in G^{\io\th}\uni \}
      \subset \wt\SX\uni, \\
\SX_{\Bm,\unip} &= \SX_{\Bm} \cap \SX\uni, \\
\pi_1^{(\Bm)} &: \wt\SX_{\Bm, \unip} \to \SX_{\Bm, \unip}, \quad (x,\Bv,gB^{\th}) \mapsto (x, \Bv).
\end{align*}
$\wt\SX_{\Bm, \unip}$ is smooth, irreducible, and the map $\pi_1^{(\Bm)}$ is 
proper, surjective.  
\par
We have the following result.

\begin{prop}  %%%%   Prop. 7.6
\begin{enumerate}
\item 
For each $\Bm \in \SQ_{n,r}$, 
\begin{equation*}
\dim \wt\CX_{\Bm,\unip} = 2n^2 - n + \sum_{i=1}^{r-1}(r-i)m_i.
\end{equation*}

\item
Assume that $\Bm \in \SQ^0_{n,r}$.  Then $\dim \wt\SX_{\Bm,\unip} = \dim \SX_{\Bm,\unip}$. 
\end{enumerate}
\end{prop}

Actually, it can be proved that if $\Bm \in \SQ^0_{n,r}$, then $\pi^{(\Bm)}_1$ is 
semi-small. 

\para{7.7.} As remarked in (6.1.1), $\SX\uni$ has infinitely many $H$-orbits if $r \ge 3$.
Thus we need to construct a certain set of subvarieties of $\SX\uni$ 
which has a similar  role as 
the set of $H$-orbits in the case where $r = 2$.  Since we want to establish the Springer 
correspondence with $W_{n,r}\wg$, these varieties must be parametrized by $\SP_{n,r}$. 
In [S3], such a variety $X_{\Bla}$ for each $\Bla \in \SP_{n,r}$ was constructed.

\begin{prop}  %%%%  Prop. 7.8
$X_{\Bla}$ is a locally closed, smooth, irreducible, $H$-stable subvariety of $\SX\uni$, 
satisfying the following properties.
\begin{enumerate}
\item 
We have
\begin{equation*}
\begin{split}
\dim X_{\Bla} = 2n^2 - 2n - &2n(\Bla) - 2n(\la^{(r-1)} + \la^{(r)})  \\
                  &+ \sum_{i=1}^{r-1}(r-i+1)|\la^{(i)}|.
\end{split}
\end{equation*}

\item
Assume that $\Bm \in \SQ_{n,r}^0$.  Then $\ol X_{\Bla(\Bm)} = \SX_{\Bm,\unip}$, 
where $\Bla(\Bm) = ((m_1), (m_2), \dots, (m_{r-1}), (m_r))$.  

\item
Asume that $\Bm \in \SQ^0_{n,r}$ and that $\Bmu \in \wt\SP(\Bm)$.  Then 
$X_{\Bmu} \subset \SX_{\Bm,\unip}$.

\item
If $r = 2$, $X_{\Bla}$ coincides with the $H$-orbit $\SO_{\Bla}$. 
\end{enumerate}
\end{prop}

\para{7.9.}
We explain the construction of $X_{\Bla}$. 
Let $\Bla \in \SP(\Bm)$ for $\Bm \in \SQ_{n,r}$.
Let $P$ be the $\th$-stable parabolic subgroup of $G$ such that 
$P^{\th}$ is the stabilizer of the (partial) isotropic flag 
$(M_{p_i})_{1 \le i \le r-2}$, and 
$L$ the $\th$-stable Levi subgroup of $P$ such that $L \supset T$. 
We shall define a set 

\begin{equation*}
\SM_{\Bla} \subset P^{\io\th}\uni \times 
           \biggl(\prod_{i=1}^{r-2}M_{p_i} \times M_{p_{r-2}}^{\perp}\biggr)
\end{equation*}
as follows; 
put $\ol M_{p_i} = M_{p_i}/M_{p_{i-1}}$ for $i = 1, \dots, r-2$, 
and $\ol M_{p_{r-1}} = M_{p_{r-2}}^{\perp}/M_{p_{r-2}}$.
Take $\Bv = (v_1, \dots, v_{r-1}) \in V^{r-1}$ such that 
$v_i \in M_{p_i}$ for $i = 1, \dots, r-2$,
and $v_{r-1} \in M_{p_{r-2}}^{\perp}$.  
Let $\ol v_i$ be the image of $v_i$ on $\ol M_{p_i}$ for $i = 1, \dots, r-1$. 
Since $\ol M_{p_{r-1}}$ has a natural symplectic structure, 
one can define an exotic symmetric space
$\SX^{(r-1)}\uni = GL(\ol M_{p_{r-1}})^{\io\th}\uni \times \ol M_{p_{r-1}}$. 
For $i = 1, \dots, r-2$, one can define an enhanced variety 
$\SX\uni^{(i)} = GL(\ol M_{p_i})\uni \times \ol M_{p_i}$, whose $GL(\ol M_{p_i})$-orbits
are parametrized by $\SP_{m_i, 2}$. 
Let $\pi_P : P^{\io\th}\uni \times V \to L^{\io\th}\uni \times V/M_{p_{r-2}}$ be 
the natural map.   
Let $\SO$ be the $H$-orbit of $z = (x, v_{r-1})$ and $\SO'$ be the $L^{\th}$-orbit
of $z' = \pi_P(z)$.  
We define $\SM_{\Bla}$ as the set of $(x, \Bv)$ satisfying the following conditions;

\begin{enumerate}

\item
$(x|_{\ol M_{p_i}}, \ol v_i) \in \SX\uni^{(i)}$ 
has type $(\la^{(i)}, \emptyset)$ for $i = 1, \dots, r-2$.

\item
$(x|_{\ol M_{p_{r-1}}}, \ol v_{r-1}) \in \SX\uni^{(r-1)}$ 
has type $(\la^{(r-1)}, \la^{(r)})$.

\item
$\SO \cap \pi_P\iv(\SO')$ is open dense in $\pi_P\iv(\SO')$.
\end{enumerate}
Finally we define 

\begin{equation*}
X_{\Bla} = \bigcup_{g \in H}g\SM_{\Bla}.
\end{equation*}

\remark{7.10.}
$\bigcup_{\Bla \in \SP_{n,r}}X_{\Bla}$ covers an open dense subset of 
$\SX\uni$, but it does not coincide with $\SX\uni$ if $r\ge 3$.
Furthermore, $X_{\Bla}$'s are not mutually disjoint in general.  

\para{7.11.}
Recall the map $\pi_1^{(\Bm)} : \wt\SX_{\Bm,\unip} \to \SX_{\Bm,\unip}$.
We define a map $\pi_{\Bm,1}: \pi_1\iv(\SX_{\Bm,\unip}) \to \SX_{\Bm, \unip}$ 
as the restriction of $\pi_1 $. 
As before, $\wt\SX_{\Bm, \unip} \subset \pi_1\iv(\SX_{\Bm, \unip})$.
Since
\begin{equation*}
(\pi_1^{(\Bm)})_*\Ql \simeq (\pi^{(\Bm)})_*\Ql|_{\SX_{\Bm,\unip}}, \qquad 
(\pi_{\Bm,1})_*\Ql \simeq (\pi_{\Bm})_*\Ql|_{\SX_{\Bm, \unip}},
\end{equation*}
$(\pi^{(\Bm)}_1)_*\Ql$ has a natural $W\nat_{\Bm}$-action inherited from 
$(\pi^{(\Bm)})_*\Ql$, and $(\pi_{\Bm,1})_*\Ql$ has a natural $W_{n,r}$-action 
inherited from $(\pi_{\Bm})_*\Ql$.  
Put $d'_{\Bm} = \dim \CX_{\Bm \unip}$. 
The following results give the Springer correspondence for $\SX\uni$, 
the former one is with respect to $W\nat_{\Bm}$, the latter one is 
with respect to $W_{n,r}$.

\begin{thm}[{[S3, Thm. 7.12, Thm. 8.17]}]  %%%%  Theorem 7.12
Assume that $\Bm \in \CQ^0_{n,r}$.  
\begin{enumerate}
\item
$(\pi_1^{(\Bm)})_*\Ql[d'_{\Bm}]$ 
is a semisimple perverse sheaf on $\CX_{\Bm,\unip}$
equipped with $W_{\Bm}\nat$-action, and is decomposed as 

\begin{equation*}
(\pi_1^{(\Bm)})_*\Ql[d'_{\Bm}] \simeq \bigoplus_{0 \le k \le m_{r-1}}
     \bigoplus_{\Bla \in \CP(\Bm(k))}
       V_{\Bla}\nat \otimes \IC(\ol X_{\Bla}, \Ql)[\dim X_{\Bla}].
\end{equation*}

\item
For $\Bla \in \SP(\Bm(k))$, we have 
\begin{equation*}
\IC(\ol X_{\Bla}, \Ql) \simeq \IC(\SX_{\Bm(k)}, \SL_{\r_{\Bla}})|_{\SX_{\Bm,\unip}}
\quad (\text{up to shift}).
\end{equation*}
\end{enumerate}
\end{thm}

\begin{thm}[{[S3, Cor. 7.4, Thm. 8.17]}] %%%%%  Theorem 7.13.
Assume that $\Bm \in \SQ^0_{n,r}$.  Then $(\pi_{\Bm,1})_*\Ql[d'_{\Bm}]$
is a semisimple perverse sheaf on $\SX_{\Bm,\unip}$ equipped with 
$W_{n,r}$-action, and is decomposed as 

\begin{equation*}
(\pi_{\Bm,1})_*\Ql[d'_{\Bm}] \simeq 
     \bigoplus_{\Bla \in \wt\SP(\Bm)}
       \wt V_{\Bla} \otimes \IC(\ol X_{\Bla}, \Ql)[\dim X_{\Bla}].
\end{equation*}
\end{thm}

\remarks{7.14.}  (i)  \ Our final goal is the Springer correspondence for $W_{n,r}$ 
given in 
Theorem 7.13.  However, the group $W_{n,r}$ and the flag variety $\SB^{\th}$ 
have no direct connection, and it is difficult to prove Theorem 7.13 directly.
On the other hand, since $W\nat_{\Bm}$ is a
prabolic subgroup of $W_n$, it has a close realtionship with $\SB^{\th}$.
Hence we first prove Theorem 7.12, and by making use of the 
relation (ii) in the theorem,  we can prove Theorem 7.13.   
\par
(ii) \ Theorem 7.12 (ii) shows that the closure $\ol X_{\Bla}$ of $X_{\Bla}$ 
is determined canonically from the decompostion of $(\pi_1^{(\Bm)})_*\Ql$. 
But there is no confidential reason for the choice of $X_{\Bla}$.  Our $X_{\Bla}$ 
is just one of such choices.

\para{7.15.}
We consider the Springer fibre for $\SX\uni$.  
For $z = (x, \Bv) \in \SX_{\Bm,\unip}$, put

\begin{align*}
\SB^{\th}_z &=  \{ gB^{\th} \in \SB^{\th} \mid g\iv xg \in B^{\io\th}, 
             g\iv \Bv \in M_n^{r-1} \}, \\
\SB_z^{\th, (\Bm)} &=  \{ gB^{\th} \in \SB^{\th} \mid g\iv xg \in B^{\io\th}, 
              g\iv \Bv \in \prod_{i=1}^{r-1}M_{p_i}  \}.
\end{align*}

Hence $\SB_z^{\th,(\Bm)} \subset \SB^{\th}_z$. 
$\SB^{\th}_z \simeq \pi_1\iv(z)$ is called the {\bf Springer fibre} of $z$, 
and $\SB_z^{\th,(\Bm)} \simeq (\pi_1^{(\Bm)})\iv(z)$ is called the 
{\bf small Springer fibre} of $z$.

\remark{7.16.}
Since $X_{\Bla}$ is not a single $H$-orbit if $r \ge 3$, it is not apriori true 
that $\dim \SB^{\th}_z, \dim \SB_z^{\th, (\Bm)}$ are constant on $X_{\Bla}$. 
In fact, this does not hold in general if $r \ge 3$. 
(But compare it with Lemma 8.9 in the enhanced case). 
\par\medskip
Assume $\Bm \in \SQ^0_{n,r}$, and $\Bla \in \wt\SP(\Bm)$.  
(Then $X_{\Bla} \subset \SX_{\Bm, \unip}$ by 7.7 (iii).) 
Put 

\begin{equation*}
d_{\Bla} = (\dim \SX_{\Bm, \unip} - \dim X_{\Bla})/2.
\end{equation*} 

\begin{lem}  %%%%  Lemma 7.17.
Assume that $\Bm \in \SQ^0_{n,r}$ and $\Bla \in \wt\SP(\Bm)$.  
For any $z \in X_{\Bla}$, we have 
\begin{equation*}
\dim \SB^{\th,(\Bm)}_z \ge d_{\Bla}.
\end{equation*} 
The set $\{ z \in X_{\Bla} \mid \dim \SB^{\th, (\Bm)}_z = d_{\Bla} \}$ forms
an open dense subset of $X_{\Bla}$. 
\end{lem}

The following result is a generalization of Corollary 5.9.

\begin{prop}  %%%%  Prop. 7.18.
Assume that $\Bm \in \SQ^0_{n,r}$, and $\Bla \in \wt\SP(\Bm)$.  
Take $z \in X_{\Bla}$
such that $\dim \SB_z^{\th, (\Bm)} = d_{\Bla}$.  Then 

\begin{enumerate}
\item
$H^{2d_{\Bla}}(\SB^{\th,(\Bm)}_z, \Ql) \simeq V\nat_{\Bla}$ as $W_{\Bm}\nat$-modules.

\item 
$H^{2d_{\Bla}}(\SB^{\th}_z, \Ql) \simeq \wt V_{\Bla}$ as $W_{n,r}$-modules.
\end{enumerate}

In particular, the map $X_{\Bla} \mapsto H^{2d_{\Bla}}(\SB^{\th}_z, \Ql)$ gives a
bijetive correspondence

\begin{equation*}
\{ X_{\Bla} \mid \Bla \in \SP_{n,r} \} \simeq W_{n,r}\wg. 
\end{equation*}
\end{prop}

\par\bigskip\medskip
\noindent
{\bf \S 8.  Enhanced variety of higher level }

\addtocounter{section}{1}
\addtocounter{thm}{-18}

\para{8.1.}
We follow the notation in \S 4.  As a generalization of the 
enhanced variety $G \times V$, we consider 
$\SX = G \times V^{r-1}$ on which $G$ acts diagonally.  
$\SX = \SX^{\text{en}}$ is called the {\bf enhanced variety of level $r$}.  
We also consider $\SX\uni = G\uni \times V^{r-1}$. 
The discussion in the case of exotic symmetric space of higher level can 
be applied also to the present situation, and in fact the arguments become simpler. 
We fix $\Bm \in \SQ_{n,r}$ (here we don't need to assume $\Bm \in \SQ^0_{n,r}$),  
and consider the varieties

\begin{align*}
\wt\SX_{\Bm} &= \{ (x, \Bv, gB) \in G \times V^{r-1} \times \SB
                     \mid g\iv xg \in B, g\iv \Bv \in \prod_{i=1}^{r-1} M_{p_i} \}, \\
   \SX_{\Bm} &= \bigcup_{g \in G}g(B \times \prod_i M_{p_i}).
\end{align*}
We define a morphism $\pi^{(\Bm)} : \wt\SX_{\Bm} \to \SX_{\Bm}$ by 
$(x,\Bv, gB) \mapsto (x, \Bv)$. 
\par
Recall that $S_{\Bm} = S_{m_1} \times \cdots \times S_{m_r}$ is a subgroup of $S_n$.
Put $d_{\Bm} = \dim \SX_{\Bm}$. 
As a generalization of (4.7.1), we have the following result.

\begin{thm}[{[S3, Thm. 4.5]}]  %%% Thm. 8.2
For each $\Bm \in \SQ_{n,r}$, $\pi^{(\Bm)}_*\Ql[d_{\Bm}]$ is a semisimple perverse 
sheaf on $\SX_{\Bm}$ equipped with $S_{\Bm}$-action, and is decomposed as 
\begin{equation*}
\pi^{(\Bm)}_*\Ql[d_{\Bm}] \simeq \bigoplus_{\r \in S_{\Bm}\wg}\r \otimes 
                            \IC(\SX_{\Bm}, \SL_{\r})[d_{\Bm}],
\end{equation*}
where $\SL_{\r}$ is a simple local system on a certain open dense subset of 
$\SX_{\Bm}$.
\end{thm}

\para{8.3.}
Put $\SX_{\Bm, \unip} = \SX_{\Bm} \cap \SX\uni$, and 
$\wt\SX_{\Bm,\unip} = (\pi^{(\Bm)})\iv(\SX_{\Bm,\unip})$.  Let 
$\pi^{(\Bm)}_1 : \wt\SX_{\Bm, \unip} \to \SX_{\Bm, \unip}$ be the restriction 
of $\pi^{(\Bm)}$ on $\wt\SX_{\Bm, \unip}$. 
We consider the complex $(\pi^{(\Bm)}_1)_*\Ql$.
As in the exotic case, $\SX\uni$ has infinitely many $G$-orbits if $r \ge 3$.
So, in order to describe the 
decomposition of $(\pi^{(\Bm)}_1)_*\Ql$, we need to introduce 
varieties $X_{\Bla}$ as in  the exotic case.    
However, in the enhanced case the situation is better than the exotic case.  
In fact we have a partition of $\SX\uni$ into pieces $X_{\Bla}$ 
indexed by $\Bla \in \SP_{n,r}$ ([S3, 5.3]), 

\begin{equation*}
\tag{8.3.1}
\SX\uni = \coprod_{\Bla \in \SP_{n,r}}X_{\Bla}
\end{equation*}
$X_{\Bla}$ is defined as follows; 
recall the notation in 4.5. Take $(x, \Bv) \in \SX\uni$ with 
$\Bv = (v_1, \dots, v_{r-1})$. Put $\ol V = V/E^xv_1$ and 
$\ol G = GL(\ol V)$.  We consider the variety 
$\SX'\uni = \ol G\uni \times \ol V^{r-2}$. Assume that 
$(x,v_1) \in G\uni \times V$ is of type $(\la^{(1)},\nu')$, where 
$\nu = \la^{(1)} + \nu'$ is the Jordan type of $x$.  
Then the type of $\ol x = x|_{\ol V}$ is $\nu'$.
Put 
$\ol\Bv = (\ol v_2, \dots, \ol v_{r-1})$, where $\ol v_i$ is the image of 
$v_i$ on $\ol V$.  Thus $(\ol x, \ol \Bv) \in \SX'\uni$.   
By induction,  we have a partition of $\SX'\uni$ 
as in (8.3.1), and there exists a unique piece $X'_{\Bla'}$ containing $(\ol x, \ol\Bv)$.
If we write $\Bla' = (\la^{(2)}, \dots, \la^{(r)})$, we obtain 
$\Bla = (\la^{(1)}, \dots, \la^{(r)}) \in \SP_{n,r}$. 
The attachment $(x,\Bv) \mapsto \Bla, \SX\uni \to \SP_{n,r}$ 
determines the set $X_{\Bla}$ as the inverse image of $\Bla$.  This gives the required 
partition. 
\par
As in the exotic case (Proposition 7.8), we have the following result.  

\begin{prop}  %%%%  Prop. 8.4.
$X_{\Bla}$ is a locally closed, smooth, irreducible, $G$-stable subvariety 
of $\SX\uni$ satisfiying the following properties.
Assume that $\Bm \in \SQ_{n,r}$.
\begin{enumerate}
\item  We have
\begin{equation*}
\dim X_{\Bla} = (n^2 - n- 2n(\Bla)) + \sum_{i=1}^{r-1}(r-i)|\la^{(i)}|.  
\end{equation*}
\item
$\SX_{\Bm,\unip} = \ol X_{\Bla(\Bm)}$, 
$($see Prop. 7.8 for the definition 
$\Bla(\Bm)$ $)$.
\item
For $\Bmu \in \SP(\Bm)$, we have $X_{\Bmu} \subset \SX_{\Bm,\unip}$.
\item
If $r = 2$, $X_{\Bla}$ coincides with the $G$-orbit $\SO_{\Bla}$. 
\end{enumerate}  
\end{prop}

\remark{8.5.} In general, the partition (8.3.1) is not compatible with 
closure relations, namely, the closure $\ol X_{\Bla}$ is not a union of 
pieces $X_{\Bmu}$.  However, we have a somewhat weaker result 
([S3, Prop. 5.11]); for each $\Bla \in \SP_{n,r}$, we have
\begin{equation*}
\ol X_{\Bla} \subset \bigcup_{\Bmu \le \Bla} X_{\Bmu}.
\end{equation*}
 
The following result gives the Springer correspondence in the case of 
the enhanced variety of higher level.  Note that (i) of the theorem 
was proved independently by Li [Li].   

\begin{thm}[{[S3, Thm. 8.13], [Li, Thm. 3.2.6]}]  %%%%  Theorem 8.6 
Assume that $\Bm \in \SQ_{n,r}$ and put $d'_{\Bm} = \dim \SX_{\Bm, \unip}$. 
\begin{enumerate}
\item
$(\pi^{(\Bm)}_1)_*\Ql[d'_{\Bm}]$ is a semisimple perverse sheaf on $\SX_{\Bm,\unip}$ 
equipped with $S_{\Bm}$-action, and is decomposed as 

\begin{equation*}
(\pi^{(\Bm)}_1)_*\Ql[d'_{\Bm}] \simeq \bigoplus_{\Bla \in \SP(\Bm)}
                \r_{\Bla}\otimes \IC(\ol X_{\Bla}, \Ql)[\dim X_{\Bla}]. 
\end{equation*}
\item 
For $\Bla \in \SP(\Bm)$, we have 
\begin{equation*}
\IC(\ol X_{\Bla}, \Ql) \simeq \IC(\SX_{\Bm}, \SL_{\r_{\Bla}})|_{\SX_{\Bm, \unip}}
\quad (\text{up to shift}). 
\end{equation*}

\end{enumerate}
\end{thm}

\para{8.7.} For $z = (x, \Bv) \in \SX_{\Bm}$, we define the Springer fibre 
$\SB^{(\Bm)}_z \simeq (\pi^{(\Bm)})\iv(z)$ by 

\begin{equation*}
\SB^{(\Bm)}_z = \{ gB \in \SB \mid g\iv xg \in B, 
             g\iv \Bv \in \prod_{i=1}^{r-1} M_{p_i} \}.
\end{equation*} 
As in the exotic case, for $\Bla \in \SP(\Bm)$ with $\Bm \in \SQ_{n,r}$, 
we define $d_{\Bla}$ by $d_{\Bla} = (\dim \SX_{\Bm, \unip} - \dim X_{\Bla})/2$. 
Then one can check easily that $d_{\Bla} = n(\Bla)$. 
By a similar argument as in the proof of Lemma 7.17, we see that 
$\dim \SB^{(\Bm)}_z \ge d_{\Bla}$.  On the other hand, since 
$\SB_z^{(\Bm)} \subset \SB_x$, and since it is well-known that 
$\dim \SB_x = n(\Bla)$, we have $\dim \SB^{(\Bm)}_z \le n(\Bla)$. 
It follows that 

\begin{lem}  %%%%  Lemma 8.8.
Assume that $\Bla \in \SP(\Bm)$.  Then for any $z \in X_{\Bla}$, 
we have $\dim \SB^{(\Bm)}_z = d_{\Bla} = n(\Bla)$. 
\end{lem}

The following result is an analogue of Propostion 7.18, and was proved 
in [S3, Prop. 8.16]. Note 
that a similar result was proved by [Li, Cor. 3.2.9] for the Borel-Moore homology. 

\begin{prop}  %%%%  Prop.  8.9.
Assume that $\Bm \in \SQ_{n,r}$, and $\Bla \in \SP(\Bm)$.  Then 
for any $z \in X_{\Bla}$, we have 
$H^{2d_{\Bla}}(\SB^{(\Bm)}_z, \Ql) \simeq \r_{\Bla}$ as $S_{\Bm}$-modules.  
\end{prop}

\para{8.10.}
We consider a generalization of the discussion on Kostka polynomials given 
in 4.10. For a general $r$, Kostka functions are given as rational functions 
$K_{\Bla, \Bmu}^{\pm}(t)$, indexed by $\Bla, \Bmu \in \SP_{n,r}$ and 
by a sign $+$ or $-$.  They are defined by fixing a total order on $\SP_{n,r}$
compatible wtih the partial order $\le$ on $\SP_{n,r}$. 
See [S1,2,4,5] for details.  The modified Kostka function 
is defined, as in the case $r = 2$, 
by $\wt K^{\pm}_{\Bla, \Bmu}(t) = t^{a(\Bmu)}K^{\pm}_{\Bla, \Bmu}(t\iv)$.  
The connection 
of those Kostka functions and the enhanced variety of higher level was discussed 
in [S4]. (Actually only the functions $K^-_{\Bla,\Bmu}(t)$ labelled by  
negative sign $-$ behave well.)  
Here we concentrate ourselves to the special case where 
$\Bmu = (-, \dots, -,\mu^{(r-1)}, \mu^{(r)})$.
Note in that case, $X_{\Bmu}$ consists of a single $G$-orbit. 
The following (partial) result can be compared to Theorem 4.11 in the $r = 2$ case.   

\begin{prop}[{[S4, Prop. 6.8]}]  %%%%  Prop. 8.11.
Take $\Bla, \Bmu \in \SP_{n,r}$, and assume that 
$\Bmu = (-,\dots, -,\mu^{(r-1)}, \mu^{(r)})$.  
Consider $ K = \IC(\ol X_{\Bla}, \Ql)$ on $\SX\uni$.  Then 
\begin{enumerate}
\item 
$\SH^i_zK = 0$ for $z \in X_{\Bmu}$ if $i$ is odd. 
\item 
Assume that $\Bmu \le \Bla$.  Then for any $z \in X_{\Bmu} \subset \ol X_{\Bla}$, 
we have

\begin{equation*}
\tag{8.11.1}
\wt K^-_{\Bla, \Bmu}(t) = t^{a(\Bla)}\sum_{i \ge 0}(\dim \SH^{\,2i}_zK)t^{ri}.
\end{equation*}
\end{enumerate}
 
\end{prop}

Note that by (8.11.1), in this special case, $\wt K^-_{\Bla,\Bmu}(t)$ and so 
$K^-_{\Bla, \Bmu}(t)$ turns out to be polynomials in $\BZ[t]$, which are independent 
of the choice of the total order. 

\para{8.12.}
An analogue of the theorem of Lascoux-Sch\"utzenberger (3.4.1) to the Kostka 
functions associated to $r$-partitions was discussed in 
[LS], [S5].  Let $\Bla = (\la^{(1)}, \dots, \la^{(r)})  \in \SP_{n,r}$.
An $r$-tuple $T = (T^{(1)}, \dots, T^{(r)})$ of tableaux is 
called a semistandard tableau of shape $\Bla$ if $T^{(i)}$ is a semistandard 
tableau of shape $\la^{(i)}$ with letters in $\{1, \dots, n\}$.
The weight $\xi = (\xi_1, \dots, \xi_n)$ of $T$ is defined by putting 
$\xi_i$ the number of letters $i$ contained in the boxes  in 
$\coprod_j T^{(j)}$ for each $i$.
For $\Bla \in \SP_{n,r}$ and $\xi \in \SP_n$, we denote by 
$SST(\Bla, \xi)$ the set of semistandard tableaux of shape $\Bla$ and weight 
$\xi$.  
Put for $\Bla \in \SP_{n,r}$,

\begin{equation*}
b(\Bla) = a(\Bla) - r\cdot n(\Bla) = |\la^{(2)}| + 2|\la^{(3)}| + \cdots + 
                   (r-1)|\la^{(r)}|. 
\end{equation*}
For $T \in SST(\Bla,\xi)$, the charge $c(T)$ is defined in 
[LS] (in the case $r = 2$), [S5] (for general $r$).  Then we have the following 
result.

\begin{thm}[{[LS, Thm. 3.12], [S5, Thm. 3.14]}] %%%% Theorem 8.14.
Let $\Bla, \Bmu \in \SP_{n,r}$, and assume that $\Bmu = (-, \dots, -, \xi)$ 
with $\xi \in \SP_n$.  Then 

\begin{equation*}
K^-_{\Bla, \Bmu}(t) = t^{b(\Bmu) - b(\Bla)}\sum_{T \in SST(\Bla, \xi)}
                               t^{r\cdot c(T)}.
\end{equation*}
\end{thm} 

%%%%
%%%%

\par\vspace{1.5cm}
\noindent
T. Shoji \\
Department of Mathematics, Tongji University \\ 
1239 Siping Road, Shanghai 200092, P. R. China  \\
E-mail: \verb|shoji@tongji.edu.cn|

\end{document}